\documentclass[11pt, a4paper, twoside,leqno]{amsart}
\usepackage[centering, totalwidth = 420pt, totalheight = 610pt]{geometry}
\usepackage{amssymb, amsmath, amsthm, enumerate, microtype, stmaryrd, url,mathrsfs}
\usepackage[latin1]{inputenc}
\usepackage[dvips, arrow, matrix, tips, curve]{xy}
\usepackage{color}
\definecolor{darkgreen}{rgb}{0,0.45,0}
\usepackage[pagebackref,colorlinks,citecolor=darkgreen,linkcolor=darkgreen]{hyperref}
\SelectTips{cm}{10}

\makeatletter
\providecommand*\xRightarrow[2][]{%
  \ext@arrow 0055{\Rightarrowfill@}{#1}{#2}}
\makeatother

\DeclareMathOperator{\ob}{ob}
\newcommand{\cat}[1]{\mathbf{#1}}
\newcommand{\op}{\mathrm{op}}
\newcommand{\id}{\mathrm{id}}
\newcommand{\thg}{{\mathord{\text{--}}}}
\newcommand{\HOM}[1]{\mathfrak{Hom}(#1)}

\newcommand{\elt}[1]{\left\llcorner{#1}\right\lrcorner}

\newcommand{\defn}[1]{\textbf{#1}}
\newcommand{\cd}[2][]{\vcenter{\hbox{\xymatrix#1{#2}}}}
\newcommand{\cdl}[2][]{\xymatrix@1#1{#2}}

\renewcommand{\b}[1]{\mathbf{#1}}
\newcommand{\A}{{\mathcal A}}
\newcommand{\B}{{\mathcal B}}
\newcommand{\C}{{\mathcal C}}
\newcommand{\D}{{\mathcal D}}

\newcommand{\R}{{\mathcal R}}
\renewcommand{\S}{{\mathcal S}}
\newcommand{\T}{{\mathcal T}}
\newcommand{\U}{{\mathcal U}}
\newcommand{\V}{{\mathcal V}}
\newcommand{\W}{{\mathcal W}}

\newcommand{\HH}{{\mathbb H}}
\newcommand{\UU}{{\mathbb U}}

\newcommand{\xtor}[1]{\cdl[@1]{{} \ar[r]|-{\object@{|}}^{#1} & {}}}
\newcommand{\tor}{\ensuremath{\relbar\joinrel\mapstochar\joinrel\rightarrow}}
\newcommand{\To}{\ensuremath{\Rightarrow}}
\newcommand{\Tor}{\ensuremath{\Relbar\joinrel\Mapstochar\joinrel\Rightarrow}}

\newcommand{\twoeq}[2][0.5]{\ar@{}[#2] \save ?(#1)*{=}\restore}
\newcommand{\twocong}[2][0.5]{\ar@{}[#2] \save ?(#1)*{\cong}\restore}
\newcommand{\rtwocell}[3][0.5]{\ar@{}[#2] \ar@{=>}?(#1)+/l 0.2cm/;?(#1)+/r 0.2cm/^{#3}}
\newcommand{\rthreecell}[3][0.5]{\ar@{}[#2] \ar@3?(#1)+/l 0.2cm/;?(#1)+/r 0.2cm/^{#3}}
\newcommand{\ltwocell}[3][0.5]{\ar@{}[#2] \ar@{=>}?(#1)+/r 0.2cm/;?(#1)+/l 0.2cm/^{#3}}
\newcommand{\ltwocello}[3][0.5]{\ar@{}[#2] \ar@{=>}?(#1)+/r 0.2cm/;?(#1)+/l 0.2cm/_{#3}}
\newcommand{\dtwocell}[3][0.5]{\ar@{}[#2] \ar@{=>}?(#1)+/u  0.2cm/;?(#1)+/d 0.2cm/^{#3}}
\newcommand{\dltwocell}[3][0.5]{\ar@{}[#2] \ar@2?(#1)+/ur  0.2cm/;?(#1)+/dl 0.2cm/^{#3}}
\newcommand{\drtwocell}[3][0.5]{\ar@{}[#2] \ar@2?(#1)+/ul  0.2cm/;?(#1)+/dr 0.2cm/^{#3}}
\newcommand{\dthreecell}[3][0.5]{\ar@{}[#2] \ar@3?(#1)+/u  0.2cm/;?(#1)+/d 0.2cm/^{#3}}
\newcommand{\drthreecell}[3][0.5]{\ar@{}[#2] \ar@3?(#1)+/ul  0.2cm/;?(#1)+/dr 0.2cm/^{#3}}
\newcommand{\dlthreecell}[3][0.5]{\ar@{}[#2] \ar@3?(#1)+/ur  0.2cm/;?(#1)+/dl 0.2cm/^{#3}}
\newcommand{\utwocell}[3][0.5]{\ar@{}[#2] \ar@{=>}?(#1)+/d 0.2cm/;?(#1)+/u 0.2cm/_{#3}}
\newcommand{\dtwocelltarg}[3][0.5]{\ar@{}#2 \ar@{=>}?(#1)+/u  0.2cm/;?(#1)+/d 0.2cm/^{#3}}
\newcommand{\utwocelltarg}[3][0.5]{\ar@{}#2 \ar@{=>}?(#1)+/d  0.2cm/;?(#1)+/u 0.2cm/_{#3}}

\newcommand{\sh}[2]{**{!/#1 -#2/}}

\newtheorem{Prop}{Proposition}

\newtheorem{Cor}[Prop]{Corollary}
\newtheorem{Thm}[Prop]{Theorem}

\theoremstyle{definition}
\newtheorem{Defn}[Prop]{Definition}
\newtheorem{Rk}[Prop]{Remark}

\renewcommand{\l}{{\mathfrak l}}
\renewcommand{\r}{{\mathfrak r}}
\renewcommand{\a}{{\mathfrak a}}

\begin{document}
\title[The low-dimensional structures formed by tricategories]{The low-dimensional structures\\formed by tricategories}
\author{Richard Garner}
\address{Department of Mathematics, Uppsala University, Box 480, S-751 06 Uppsala,
       Sweden}\email{richard@math.uu.se}
       \author{Nick Gurski}
       \address{Department of Mathematics, Yale University, 10 Hillhouse Avenue,
        PO Box 208283, New Haven, CT 06520-8283, USA}
       \email{michaeln.gurski@yale.edu}
\thanks{The first-named author acknowledges the support of a Marie Curie Intra-European
Fellowship, Project No.\ 040802, and a Research Fellowship of St John's
College, Cambridge.}
    \maketitle

\begin{abstract}
We form tricategories and the homomorphisms between them into a bicategory,
whose $2$-cells are certain degenerate tritransformations. We then enrich this
bicategory into an example of a three-dimensional structure called a
\emph{locally cubical bicategory}, this being a bicategory enriched in the
monoidal 2-category of pseudo double categories. Finally, we show that every
sufficiently well-behaved locally cubical bicategory gives rise to a
tricategory, and thereby deduce the existence of a tricategory of
tricategories.
\end{abstract}

\section{Introduction}
A major impetus behind many developments in 2-dimensional category theory has
been the observation that, just as the fundamental concepts of set theory are
categorical in nature, so the fundamental concepts of category theory are
2-categorical in nature. In other words, if one wishes to study categories ``in
the small''~--~as mathematical entities in their own right rather than as
universes of discourse~--~then a profitable way of doing this is by studying
the 2-categorical properties of $\cat{Cat}$, the 2-category of all
categories.\footnote{Here, and elsewhere, we will adopt a common-sense attitude
to set-theoretic issues, assuming a sufficient supply of Grothendieck universes
and leaving it to the reader to qualify entities with suitable constraints on
their size.}

Once one moves from the study of categories to the study of (possibly weak)
$n$-categories, it is very natural to generalise the above maxim, and to assert
that \emph{the fundamental concepts of $n$-category theory are
$(n+1)$-categorical in nature.} This is a profitable thing to do: for example,
consider the coherence theorem for bicategories \cite{MLP}, which in its
simplest form states that
\[\text{\emph{Every bicategory is biequivalent to a 2-category.}}\]
\textit{A priori}, this is merely a statement about individual bicategories;
but we may also read it as a statement about the tricategory of bicategories
$\cat{Bicat}$, since ``biequivalent'' may be read as ``internally biequivalent
in the tricategory $\cat{Bicat}$.''\footnote{In practice, one would tend to use
the local definition of biequivalence, wherein $\B$ is biequivalent to $\B'$ if
there exists a homomorphism $F \colon \B \to \B'$ which is biessentially
surjective on objects and locally an equivalence of categories; but as long as
we assume the axiom of choice, the difference between the two definitions is
merely one of presentation.}\ \  Thus another way of stating the above would be
to say that the 2-categories are biequivalence-dense in $\cat{Bicat}$.

This maxim permeates almost all research in higher-dimensional category theory,
and so we draw attention to it here, not in order to point out where we might use
it, but rather where we might \emph{not} use it. For instance, consider once
again the coherence theorem for bicategories. We may restate it slightly more
tightly as:
\begin{gather*}\text{\emph{Every bicategory is biequivalent to a
2-category}}\\\text{\emph{via an identity-on-objects
biequivalence.}}\end{gather*} The restriction to identity-on-objects
biequivalences affords us an interesting simplification, since, as pointed out
in \cite{LP}, we can express such a biequivalence as a mere \emph{equivalence}
in a suitable 2-category, which we denote by $\cat{Bicat}_2$. The 0-cells of
$\cat{Bicat}_2$ are the bicategories; the 1-cells are the homomorphisms between
them; and the 2-cells are the \emph{icons} of~\cite{Lackicons}. These are
degenerate oplax natural transformations whose every $1$-cell component is an
identity: we will meet them in more detail in Section~\ref{bicat-tricat} below.

With the help of the 2-category $\cat{Bicat}_2$, the coherence theorem for
bicategories can be made into a 2-categorical, rather than a tricategorical,
statement: namely, that the 2-categories are equivalence-dense in
$\cat{Bicat}_2$ (cf. Theorem 5.4 of \cite{LP}). This is a somewhat tighter
result; moreover, the 2-category $\cat{Bicat}_2$ is much simpler to work with
than the tricategory $\cat{Bicat}$. Thus we should revise our general maxim,
and acknowledge that \emph{some} of the fundamental concepts of $n$-category
theory may be expressible using fewer than $(n+1)$ dimensions. Consequently,
when we study $n$-categories, it may be useful to form them not only into an
$(n+1)$-category, but also into suitable lower-dimensional structures. It is
the purpose of this paper to do this in the case $n=3$. We construct both a
bicategory of tricategories $\cat{Tricat}_2$ and a tricategory of tricategories
$\cat{Tricat}_3$: where in both cases, the 2-cells are suitably scaled-up
analogues of the bicategorical icons mentioned above.

In \cite{Lackicons}, Lack gives a number of motivations for studying the
2-category $\cat{Bicat}_{2}$ of bicategories, lax functors, and icons.  Many of
these motivations have obvious analogues one dimension higher.  For instance,
the coherence theorem for tricategories can be restated as
\[\text{\emph{Every tricategory is internally biequivalent to a $\cat{Gray}$-category
in the tricategory $\cat{Tricat}_{3}$.}}\] On the other hand, coherence for
tricategories internal to the \textit{bicategory} $\cat{Tricat}_{2}$ is an open
question.  The structures $\cat{Tricat}_{2}$ and $\cat{Tricat}_{3}$ also
provide avenues for studying the simplicial nerves of tricategories, thus
allowing comparisons with work of Street \cite{streetomega} to be pursued in
dimension three.  Moreover, it is shown in~\cite{Lackicons} that the
$2$-category of monoidal categories embeds nicely in $\cat{Bicat}_{2}$; and
similarly, we show that the tricategory of monoidal bicategories~--~as
constructed in \cite{cheng-gurski}~--~embeds nicely in $\cat{Tricat}_{3}$.

An outline of the paper is as follows.  In Section~\ref{bicat-tricat}, we
construct a bicategory of tricategories, $\cat{Tricat}_{2}$.  The construction
is straightforward and computational. The 2-cells of this bicategory we call
\textit{ico-icons}: they can be seen as doubly degenerate oplax
tritransformations whose $0$- and $1$-cell components are identities. More
explicitly, they exist only between trihomomorphisms which agree on $0$- and
$1$-cells, and are given by a collection of (not necessarily invertible)
$3$-cell components together with coherence data and axioms.

In Section~\ref{Sec:towards}, we describe a tricategory of tricategories,
$\cat{Tricat}_3$. The first candidate we consider for its $2$-cells are the
\emph{oplax icons}, which are singly degenerate oplax tritransformations: they
exist only between trihomomorphisms which agree on $0$-cells, and are given by
a collection of (not necessarily invertible) $2$- and $3$-cell components
together with coherence data and axioms. These generalise the $2$-cells of
$\cat{Tricat}_2$, since every ico-icon is an oplax icon: indeed, the ico-icons
are precisely the \emph{i}dentity \emph{c}omponents \emph{o}plax icons.
However, oplax icons turn out to be too lax to compose properly: the same
phenomenon which occurs if one tries to replace the weak transformations in the
tricategory $\cat{Bicat}$ with oplax transformations. Thus instead we take the
$2$-cells of $\cat{Tricat}_3$ to be the smaller class of \emph{pseudo-icons}:
these being oplax icons whose $3$-dimensional data is invertible.

Although we describe the tricategory $\cat{Tricat}_3$ in
Section~\ref{Sec:towards}, we do not complete its construction. One reason is
that we want to avoid giving unenlightening tricategorical coherence
computations as far as possible, to which end, we would like to reuse the work
we did in Section~\ref{bicat-tricat}; and though intuitively this is not a
problem, technically it is rather troublesome. A second reason is that we wish
to explain an unusual discrepancy, namely that the bicategory $\cat{Tricat}_2$
carries some information which the tricategory $\cat{Tricat}_3$ cannot, in that
an ico-icon ($2$-cell of $\cat{Tricat}_2$) cannot be viewed as a pseudo-icon
($2$-cell of $\cat{Tricat}_3$) unless it is \emph{invertible}.

In Sections~\ref{Sec:locallydouble}--\ref{Sec:fromltt} we describe a general
mechanism which allows us to clear up both of the above issues. This begins in
Section~\ref{Sec:locallydouble} with the introduction of a new kind of
three-dimensional categorical structure which we call a \emph{locally cubical
bicategory}. Like a tricategory, it has $0$-, $1$-, $2$-, and $3$-cells; but
the $2$-cells come in two different kinds, vertical and horizontal, whilst the
$3$-cells are cubical in nature. Moreover, the coherence axioms that are to be
satisfied are of a bicategorical, rather than a tricategorical kind, and so the
resultant structure is computationally more tractable than a tricategory. As a
first application of this theory, we are able to show quite easily that the
totality of bicategories (and more generally the totality of \emph{pseudo
double categories} in the sense of~\cite{GP1}) form a locally cubical
bicategory.
%

Section~\ref{Sec:local} then describes a locally cubical bicategory of
tricategories which we denote by $\mathfrak{Tricat}_3$. The construction is
once again straightforward and computational, and reuses the work done in
Section~\ref{bicat-tricat}. The objects and $1$-cells of $\mathfrak{Tricat}_3$
are just tricategories and trihomomorphisms; the vertical $2$-cells are the
ico-icons from $\cat{Tricat}_2$; the horizontal $2$-cells are the pseudo-icons
from $\cat{Tricat}_3$; whilst the $3$-cells are ``cubical icon modifications''.
In particular, $\mathfrak{Tricat}_3$ is a rich enough structure to encode all
the information from both $\cat{Tricat}_2$ and $\cat{Tricat}_3$. This resolves
the second of the issues mentioned above.

In order to resolve the first issue, we appeal to a general theory which allows
us to construct tricategories out of sufficiently well-behaved locally cubical
bicategories: more precisely, those with the property that every invertible
vertical $2$-cell gives rise to a horizontal $2$-cell. This general theory is
described in detail in Section~\ref{Sec:fromltt}; whilst in
Section~\ref{tritri}, we are able to apply it to the locally cubical bicategory
$\mathfrak{Tricat}_3$, thereby deducing the existence of the tricategory of
tricategories $\cat{Tricat}_3$. Additionally, we identify the tricategory of
monoidal bicategories inside of $\cat{Tricat}_{3}$.

\textbf{Notation.}\ \ \ We follow \cite{bicats} and \cite{review} where it
concerns 2-\ and bicategories: so in particular, our oplax natural
transformations $\alpha \colon F \Rightarrow G$ have 2-cell components given by
$\alpha_f \colon \alpha_B . Ff \Rightarrow Gf . \alpha_A$. We will tend to use
either juxtaposition or the connective ``$.$'' to denote composition, relying
on context to sort out precisely which sort of composition is intended. When it
comes to tricategories, our primary references are \cite{GPS} and
\cite{nicktricats}, but with a preference for the ``algebraic'' presentation of
the latter: though we will not use this algebraicity in any essential way.

We will also make use of \emph{pasting diagrams} of 2-cells inside
bicategories. Such diagrams are only well-defined up to a choice of bracketing
of their boundary, and so we assume such a choice to have been made wherever
necessary. Occasionally we will need to use similar pasting diagrams of 2-cells
in a tricategory, and the same caveat holds, only more so: here, the diagram is
only well-defined up to a choice of order in which the pasting should be
performed; and again, we assume such a choice to have been made. We adopt one
further convention regarding pasting diagrams. Suppose we are given a 2-cell
$\alpha \colon h(gf) \Rightarrow h'(g'f')$ in a bicategory $\B$, thus:
\[
\cd{
  & X \ar[r]^g & Y \ar[dr]^h \\
W \ar[ur]^f \ar[dr]_{f'} \dtwocell{rrr}{\alpha} & & & Z\text, \\
  & X' \ar[r]_{g'} & Y' \ar[ur]_{h'}}
\]
together with a homomorphism of bicategories $F \colon \B \to \C$. Applying $F$
to $\alpha$ yields a 2-cell $F(h(gf)) \Rightarrow F(h'(g'f'))$ of $\C$, but
frequently, we will be more interested in the 2-cell
\[
\cd{
  & FX \ar[r]^{Fg} & FY \ar[dr]^{Fh} \\
FW \ar[ur]^{Ff} \ar[dr]_{Ff'} \dtwocell{rrr}{} & & & FZ\text; \\
  & FX' \ar[r]_{Fg'} & FY' \ar[ur]_{Fh'}}
\]
obtained by pasting $F\alpha$ with suitable coherence constraints for the
homomorphism $F$: and we will consistently denote the 2-cell obtained in this
way by $\overline{F\alpha}$.

\section{A bicategory of tricategories}\label{bicat-tricat}

We begin by describing the lowest-dimensional structure into which
tricategories and their homomorphisms can form themselves. At first, one might
think that this would be a category; but unfortunately, composition of
trihomomorphisms fails to be associative on the nose, as it requires one to
compose 1-cells in a hom-bicategory, which is itself not an associative
operation. Consequently, the best we can hope for is a \emph{bicategory} of
tricategories, which we will denote by $\cat{Tricat}_2$.

The simplest such bicategory would have trihomomorphisms as its 1-cells and
\emph{blips} as its 2-cells. According to \cite{GPS}, blips are very degenerate
tritransformations which can only exist between two trihomomorphisms $F, G
\colon \S \to \T$ which agree on 0-, 1-,~2-\ and 3-cells. Though one might
think that this forces $F$ and $G$ to be the same, they can in fact differ with
respect to certain pieces of coherence data: and a ``blip'' is the means by
which one measures these differences.

However, if we are going to form a bicategory of tricategories, it may as well
be the most general possible one; and so we will consider more general sorts of
both 1-~and 2-cells. Let us begin by looking at the 1-cells.
\begin{Defn}\label{laxhomom}
Let $\S$ and $\T$ be tricategories. A \defn{lax homomorphism} $F \colon \S \to
\T$ is a lax morphism of tricategories in the sense of \cite{GPS}, all of whose
coherence 3-cells are invertible. Hence $F$ consists of:
\begin{itemize}
\item A function $F \colon \ob \S \to \ob \T$;
\item Homomorphisms of bicategories $F_{A, B} \colon \S(A, B) \to \T(FA,
    FB)$;
\item 2-cells $\iota_A \colon I_{FA} \to FI_A$;
\item 2-cells $\chi_{f, g} \colon Fg.Ff \Rightarrow F(gf)$, pseudo-natural
    in $f$ and $g$;
\item Invertible modifications $\omega$, $\delta$ and $\gamma$ witnessing
    the coherence of $\iota$ and $\chi$,
\end{itemize}
all subject to the axioms for a morphism of tricategories as found in
\cite{GPS}.
\end{Defn}

The notion of lax homomorphism is a sensible one from many angles.  We can
compose lax homomorphisms just as we would compose homomorphisms of
tricategories. If we are given a pair of monoidal bicategories
\cite{monoidalbicats} which we view as one-object tricategories, then the lax
homomorphisms between them are the natural bicategorical generalisation of a
lax monoidal functor (\emph{weak monoidal homomorphisms}, in the
terminology of \cite{monoidalbicats}).  Lax homomorphisms from the terminal
tricategory into $\T$ classify \emph{pseudomonads} in $\T$~--~that is, monads
whose associativity and unit laws have been weakened to hold up to coherent
isomorphism, and in a similar vein we may use lax homomorphisms to give a
succinct definition of an \emph{enriched bicategory} in the sense of
\cite{Carmody, Lack}~--~that is, of a bicategory ``enriched in a tricategory'',
which is a one-dimension-higher version of a category enriched in a
bicategory~\cite[Section 5.5]{bicats}, which is in turn a generalisation of the
familiar notion of a category enriched in a monoidal category. We shall see a
little more of enriched bicategories in Section~\ref{Sec:locallydouble}.

We now turn to the 2-cells of $\cat{Tricat}_2$. The most informative precedent
here is the corresponding notion one dimension down: the \emph{icons} of
\cite{Lackicons, LP}. As mentioned in the Introduction, these are degenerate
oplax transformations between homomorphisms of bicategories which agree on
0-cells. To be precise, given two such homomorphisms of bicategories $F, G
\colon \B \to \C$, an \defn{icon} $\alpha \colon F \To G$ is given by
specifying for each 1-cell $f \colon A \to B$ of $\B$, a 2-cell $\alpha_f
\colon Ff \To Gf$ of $\C$ such that:
\begin{itemize}
\item For each 2-cell $\sigma \colon f \To g$ of $\B$, the following
    diagram commutes:
\[\cd{
 Ff \ar@2[r]^{\alpha_f} \ar@2[d]_{F \sigma} & Gf \ar@2[d]^{G \sigma} \\
 Fg \ar@2[r]_{\alpha_g} & Gg\text;
}\]
\item For each object $A \in \B$, the following diagram commutes:
\[\cd{
 {\id_{FA}} \ar@2[r]^{\cong} \ar@{=}[d] & F\id_A \ar@2[d]^{\alpha_{\id_A}} \\
 {\id_{GA}} \ar@2[r]_{\cong} & G\id_A\text;
}\]
\item For each pair of composable 1-cells $f \colon A \to B$, $g \colon B
    \to C$ in $\B$, the following diagram commutes:
\[\cd{
 Fg . Ff \ar@2[r]^{\cong} \ar@2[d]_{\alpha_{g} . \alpha_{f}} & F(gf) \ar@2[d]^{\alpha_{gf}} \\
 Gg . Gf \ar@2[r]_{\cong} & G(gf)\text,
}\]
\end{itemize}
where the arrows labelled with $\cong$ witness the pseudo-functoriality
of $F$ and $G$. There is a bijection between icons $F \Rightarrow G$ and those
oplax natural transformations $F \Rightarrow G$ whose components are all
identities (hence the name: \emph{i}dentity \emph{c}omponent \emph{o}plax
\emph{n}atural transformation); however, icons differ crucially from the oplax
natural transformations representing them in regard to the manner of their
\emph{composition}. Indeed, composition of oplax natural transformations is
only associative up to invertible modification, whilst icons admit a strictly
associative composition; and it is this which allows bicategories,
homomorphisms and icons to form a $2$-category $\cat{Bicat}_2$.

The 2-cells of $\cat{Tricat}_2$ we are about to describe~--~the
\emph{ico-icons}~--~can be seen as higher-dimensional analogues of these
bicategorical icons. They are doubly degenerate oplax tritransformations
between lax trihomomorphisms which agree on both $0$- and $1$-cells. Here
again, composition of ico-icons will not simply be composition of
tritransformations, but rather a modified form of that composition which is
\emph{strictly} associative. The choice of the name ico-icon will be explained
by Proposition~\ref{embed} below.

\begin{Defn}\label{icoicon}
Given lax homomorphisms $F$, $G \colon \S \to \T$, an \defn{ico-icon} $\alpha
\colon F \Rightarrow G$ may exist only if $F$ and $G$ agree on objects and
1-cells of $\S$; and is then given by the following data:
\begin{enumerate}[(TD1)]
\item[(TD1)] For each pair of objects $A, B \in \S$, an icon
\[\alpha_{A, B} \colon F_{A, B} \Rightarrow G_{A, B} \colon \S(A, B) \to \T(FA, FB)\]
(so in particular, for each 2-cell $\theta \colon f \Rightarrow g$ of $\S$,
a 3-cell of $\T$:
\[\cd{
    Ff \ar@2[r]^{F \theta} \ar@{=}[d] \dthreecell{dr}{\alpha_{\theta}}& Fg \ar@{=}[d] \\
    Gf \ar@2[r]_{G \theta} & Gg
}\quad\text{);}\]
\item[(TD2)] For each object $A$ of $\S$, a 3-cell of $\T$:
\[
\cd{
    I_{FA} \ar@2[r]^{\iota^F_A} \ar@{=}[d] \dthreecell{dr}{M^{\alpha}_A} & FI_A \ar@{=}[d] \\
    I_{GA} \ar@2[r]_{\iota^G_A} & GI_A\text;
}\]
\item[(TD3)] For each pair of composable 1-cells $f \colon A \to B$, $g
    \colon B \to C$ of $\S$, a 3-cell of $\T$:
\[
\cd{
    Fg . Ff \ar@2[r]^{\chi^F_{f, g}} \ar@{=}[d] \dthreecell{dr}{\Pi^{\alpha}_{f, g}} & F(g f) \ar@{=}[d] \\
    Gg . Gf \ar@2[r]_{\chi^G_{f, g}} & G(g f)\text;
}\]
\end{enumerate}
subject to the following axioms:
\begin{enumerate}[(TA1)]
\item[(TA1)] For each pair of 2-cells $\theta \colon f \Rightarrow g \colon
    B \to C$ and $\theta' \colon f' \Rightarrow g' \colon A \to B$ of $\S$,
    the following pasting equality holds:
\[ \hphantom{=} \quad \cd[@!@-0.5em]{
    & Ff . Ff'
      \ar@2[r]^{\chi^F} \ar@{=}[dl] \dlthreecell{d}{\Pi^{\alpha}}
    & F(ff')
      \ar@{=}[dl] \ar@2[dr]^{F(\theta \theta')} \dlthreecell{dd}{\alpha_{\theta \theta'}}
    \\ Gf . Gf'
      \ar@2[r]^{\chi^G} \ar@2[dr]_{G\theta . G\theta'}
    & G(ff')
      \ar@2[dr]^{G(\theta \theta')} \drthreecell{d}{\chi^G}
    & & F(g g')
      \ar@{=}[dl] \\
    & Gg . Gg'
      \ar@2[r]_{\chi^G}
    & G(gg')
}\]\[ = \quad \cd[@!@-0.5em]{
    & Ff . Ff' \ar@2[r]^{\chi^F} \ar@{=}[dl] \ar@2[dr]_{F\theta . F\theta'}
      \dlthreecell{dd}{\alpha_\theta . \alpha_{\theta'}}
    & F(ff') \ar@2[dr]^{F(\theta \theta')} \drthreecell{d}{\chi^F}
    \\ Gf . Gf'
      \ar@2[dr]_{G\theta . G\theta'} & & Fg . Fg'
      \ar@2[r]^{\chi^F} \ar@{=}[dl] \dlthreecell{d}{\Pi^\alpha} &
    F(gg')\text;
      \ar@{=}[dl] \\
    & Gg . Gg' \ar@2[r]_{\chi^G} & G(gg')\text. }\]
\item[(TA2)] For each 1-cell $f \colon A \to B$ of $\S$, the following
    pasting equality holds:
\[
 \hphantom{=} \quad \cd[@!C@!R@R+1em]{
    & FI_{B}.Ff
      \ar@2[r]^{\chi^F} \dthreecell{ddr}{\gamma^F}
    & F(I_B.f)
      \ar@2[dr]^{F\l}
    \\ I_{FB}.Ff
      \ar@2[dr]^{\l} \ar@2[ur]^{\iota^F.1} \ar@{=}[d]
    & & & Ff
      \ar@{=}[d]
    \\ I_{GB}.Gf
      \ar@2[dr]_{\l} \twoeq{r}
    & Ff
      \ar@{=}[r] \ar@{=}[d] \twoeq{dr}
    & Ff
      \ar@{=}[ur] \ar@{=}[d] \twoeq{r}
    & Gf
    \\ & Gf
      \ar@{=}[r]
    & Gf
      \ar@{=}[ur]
}\]
\[
= \quad \cd[@!C@!R@R+1em]{
    & FI_B . Ff
      \ar@2[r]^{\chi^F} \ar@{=}[d] \dthreecell{dr}{\Pi^\alpha}
    & F(I_B . f)
      \ar@2[dr]^{F\l} \ar@{=}[d]
    \\ I_{FB}.Ff
      \ar@2[ur]^{\iota^F.1} \ar@{=}[d] \ar@{}[r] \ar@3?(0.5)+/u  0.2cm/;?(0.5)+/d 0.2cm/_(1.5)*+[l]{\scriptstyle \ \ M^\alpha.1}
    & GI_B . Gf
      \ar@2[r]_{\chi^G} \dthreecell{ddr}{\gamma^G}
    & G(I_B . f)
      \ar@2[dr]_{G \l} \dthreecell{r}{\alpha_\l}
    & Ff
      \ar@{=}[d]
    \\ I_{GB}.Gf
      \ar@2[dr]_{\l} \ar@2[ur]_{\iota^G.1}
    & & & Gf\text.
    \\ & Gf
      \ar@{=}[r]
    & Gf
      \ar@{=}[ur]
}\]
\item[(TA3)] For each 1-cell $f \colon A \to B$ of $\S$, the following
    pasting equality holds:
\[
 \hphantom{=} \quad \cd[@!C@!R@R+1em]{
    & Ff.FI_{A}
      \ar@2[r]^{\chi^F} \dthreecell{ddr}{\delta^F}
    & F(f.I_A)
      \ar@2[dr]^{F\r}
    \\ Ff.I_{FA}
      \ar@2[dr]^{\r} \ar@2[ur]^{1.\iota^F} \ar@{=}[d]
    & & & Ff
      \ar@{=}[d]
    \\ Gf.I_{GA}
      \ar@2[dr]_{\r} \twoeq{r}
    & Ff
      \ar@{=}[r] \ar@{=}[d] \twoeq{dr}
    & Ff
      \ar@{=}[ur] \ar@{=}[d] \twoeq{r}
    & Gf
    \\ & Gf
      \ar@{=}[r]
    & Gf
      \ar@{=}[ur]
}\]
\[
= \quad \cd[@!C@!R@R+1em]{
    & Ff.FI_A
      \ar@2[r]^{\chi^F} \ar@{=}[d] \dthreecell{dr}{\Pi^\alpha}
    & F(f.I_A)
      \ar@2[dr]^{F\r} \ar@{=}[d]
    \\ Ff.I_{FA}
      \ar@2[ur]^{1.\iota^F} \ar@{=}[d] \ar@{}[r] \ar@3?(0.5)+/u  0.2cm/;?(0.5)+/d 0.2cm/_(1.5)*+[l]{\scriptstyle \ \ 1.M^\alpha}
    & Gf.GI_A
      \ar@2[r]_{\chi^G} \dthreecell{ddr}{\delta^G}
    & G(f.I_A)
      \ar@2[dr]_{G \r} \dthreecell{r}{\alpha_\r}
    & Ff
      \ar@{=}[d]
    \\ Gf.I_{GA}
      \ar@2[dr]_{\r} \ar@2[ur]_{1.\iota^G}
    & & & Gf\text.
    \\ & Gf
      \ar@{=}[r]
    & Gf
      \ar@{=}[ur]
}\]
\item[(TA4)] For each triple $f, g, h$ of composable 1-cells of $\S$, the
    following pasting equality holds:
\[
 \hphantom{=} \quad \cd[@!C@!R@R+1em]{
    & F(hg). Ff
      \ar@2[r]^{\chi^F} \dthreecell{ddr}{\omega^F}
    & F((hg)f)
      \ar@2[dr]^{F\a}
    \\ (Fh.Fg).Ff
      \ar@2[dr]^{\a} \ar@2[ur]^{\chi^F.1} \ar@{=}[d]
    & & & F(h(gf))
      \ar@{=}[d]
    \\ (Gh.Gg).Gf
      \ar@2[dr]_{\a} \twoeq{r}
    & Fh.(Fg. Ff)
      \ar@2[r]^{1. \chi^F} \ar@{=}[d] \dthreecell{dr}{1.\Pi^\alpha}
    & Fh. F(gf)
      \ar@2[ur]^{\chi^F} \ar@{=}[d] \dthreecell{r}{\Pi^\alpha}
    & G(h(gf))
    \\ & Gh.(Gg.Gf)
      \ar@2[r]_{1.\chi^G}
    & Gh.G(gf)
      \ar@2[ur]_{\chi^G}
}\]
\[
= \quad \cd[@!C@!R@R+1em]{
    & F(hg). Ff
      \ar@2[r]^{\chi^F} \ar@{=}[d] \dthreecell{dr}{\Pi^\alpha}
    & F((hg)f)
      \ar@2[dr]^{F\a} \ar@{=}[d]
    \\ (Fh.Fg).Ff
      \ar@2[ur]^{\chi^F.1} \ar@{=}[d] \ar@{}[r] \ar@3?(0.5)+/u  0.2cm/;?(0.5)+/d 0.2cm/_(1.5)*+[l]{\scriptstyle \ \ \Pi^\alpha.1}
    & G(hg).Gf
      \ar@2[r]_{\chi^G} \dthreecell{ddr}{\omega^G}
    & G((hg)f)
      \ar@2[dr]_{G \a} \dthreecell{r}{\alpha_\a}
    & F(h(gf))
      \ar@{=}[d]
    \\ (Gh.Gg).Gf
      \ar@2[dr]_{\a} \ar@2[ur]_{\chi^G.1}
    & & & G(h(gf))\text.
    \\ & Gh.(Gg.Gf)
      \ar@2[r]_{1.\chi^G}
    & Gh.G(gf)
      \ar@2[ur]_{\chi^G}
}\]
\end{enumerate}
\end{Defn}

Observe that, because the raw data for an ico-icon is a collection of 3-cells
in the target tricategory, there is no possibility of introducing a third
dimension of structure given by ``ico-icon modifications''. To do this we have
to look at \emph{singly} degenerate, rather than \emph{doubly} degenerate,
oplax tritransformations. We do this in the next Section.

Now, in order to show that this collection of 0-, 1-~and 2-cells forms a
bicategory, we have to give additional \emph{data}~--~vertical composition of
2-cells, horizontal composition of 1-\ and 2-cells and associativity and
unitality constraints~--~subject to additional \emph{axioms}~--~the category
axioms for vertical composition, the middle-four interchange axiom and the
pentagon and triangle axioms for the associativity and unit constraints.

We start with the vertical structure: the identity 2-cell $\id_F \colon F
\Rightarrow F$ in $\cat{Tricat}_2$ we take to be given by the following data:
\[(\id_F)_{A, B} = \id_{F_{A, B}}\text, \quad
M^{\id_F}_A = \id_{\iota^F_A} \quad \text{and} \quad \Pi^{\id_F}_{f, g} = \id_{\chi^F_{f, g}}
\text.\] Each of the axioms (TA1)--(TA4) now expresses that something is equal to itself pasted
together with some identity 3-cells, which is clear enough. Next, given 2-cells
$\alpha \colon F \Rightarrow G$ and $\beta \colon G \Rightarrow H$ in
$\cat{Tricat}_2$, we take $\beta\alpha \colon F \Rightarrow H$ to be given by
the following data:
\[(\beta \alpha)_{A, B} = \beta_{A, B} . \alpha_{A, B}\text, \quad
M^{\beta\alpha}_A = M^{\beta}_A . M^{\alpha}_A \quad \text{and} \quad \Pi^{\beta\alpha}_{f, g} =
\Pi^{\beta}_{f, g} . \Pi^{\alpha}_{f, g}\text.\] Each of the axioms (TA1)--(TA4) for this data
follow from juxtaposing the corresponding axioms for $\alpha$ and $\beta$ in a
very straightforward manner. Moreover, because vertical composition of 3-cells
in a tricategory is strictly associative and unital, so is the vertical
composition of 2-cells in $\cat{Tricat}_2$.

We turn now to the horizontal structure. Horizontal identities and composition
for 1-cells are the identities and composition for lax homomorphisms as
detailed in \cite{Lack}; whilst given 2-cells $\alpha \colon F \Rightarrow F'
\colon \S \to \T$ and $\beta \colon G \Rightarrow G' \colon \T \to \U$, their
horizontal composite $\beta \ast \alpha \colon GF \Rightarrow G'F' \colon \S
\to \U$ is given by:
\begin{enumerate}[(TD1)]
\item[(TD1)] $(\beta \ast \alpha)_{A, B} := \beta_{A, B} \ast \alpha_{A,
    B}$, where $\ast$ on the right-hand side is the horizontal composite of
    the underlying icons in the 2-category $\cat{Bicat}_2$ of the
    Introduction. In particular, given a 2-cell $\theta \colon f
    \Rightarrow g$ of $\S$, we have
\[(\beta \ast \alpha)_\theta \quad = \quad \cd{
    GFf \ar@2[r]^{GF \theta} \ar@{=}[d] \dthreecell{dr}{\beta_{F\theta}}& GFg \ar@{=}[d] \\
    G'Ff \ar@2[r]^{G'F \theta} \ar@{=}[d] \dthreecell{dr}{G'\alpha_{\theta}}& G'Fg \ar@{=}[d] \\
    G'F'f \ar@2[r]_{G'F' \theta} & G'F'g\text{;}
}\]
\item[(TD2)]
\[
M^{\beta \ast \alpha}_A \quad := \quad \cd{
    I_{GFA}
      \ar@2[r]^{\iota^G} \ar@{=}[d] \dthreecell{dr}{M^{\beta}_{FA}}
    & GI_{FA}
      \ar@{=}[d] \ar@2[r]^{G\iota^F} \dthreecell{dr}{\beta_{\iota^F}}
    & GFI_A
      \ar@{=}[d]
    \\ I_{G'FA}
      \ar@2[r]^{\iota^{G'}} \ar@{=}[d] \twoeq{dr}
    & G'I_{FA}
      \ar@{=}[d] \ar@2[r]^{G'\iota^F} \dthreecell{dr}{G'M^{\alpha}_A}
    & G'FI_A
      \ar@{=}[d]
    \\ I_{G'F'A}
      \ar@2[r]_{\iota^{G'}}
    & G'I_{F'A}
      \ar@2[r]_{G'\iota^{F'}}
    & G'F'I_A\text,}
\]
\item[(TD3)]
\[
\Pi^{\beta \ast \alpha}_{f, g} \quad := \quad \cd{
    GFg . GFf
      \ar@2[r]^{\chi^G} \ar@{=}[d] \dthreecell{dr}{\Pi^{\beta}_{Ff, Fg}}
    & G(Fg . Ff)
      \ar@2[r]^{G\chi^F} \ar@{=}[d] \dthreecell{dr}{\beta_{\chi^F}}
    & GF(gf)
      \ar@{=}[d] \\
    G'Fg . G'Ff
      \ar@2[r]^{\chi^{G'}} \ar@{=}[d] \twoeq{dr}
    & G'(Fg . Ff)
      \ar@2[r]^{G'\chi^F} \ar@{=}[d] \dthreecell{dr}{G'\Pi^\alpha_{f, g}}
    & G'F(gf)
      \ar@{=}[d] \\
    G'F'g . G'F'f
      \ar@2[r]_{\chi^{G'}}
    & G'(F'g . F'f)
      \ar@2[r]_{G'\chi^{F'}}
    & G'F'(gf)\text.
}\]
\end{enumerate}
We must check that these data satisfy (TA1)--(TA4)\label{checkaxioms}. If we
view the pasting equalities in these axioms as equating two ways round a cube
or a hexagonal prism, then this verification is a matter of taking a suitable
collection of such cubes and prisms for $\beta$ and $\alpha$ and sticking them
together in the right way. When realised in two dimensions, this amounts to
displaying a succession of equalities of rather large pasting diagrams. We
leave the task of reconstructing these to the reader.

Let us consider now the middle-four interchange axiom. Asking for this be
satisfied amounts to checking that the other obvious way of defining $\beta
\ast \alpha$~--~via $GF'$ rather than $G'F$~--~gives the same answer; and this
follows quickly from the middle-four interchange law in the hom-bicategories of
$\U$, and the first icon axiom for $\beta$.

It remains to give the associativity and unit constraints $a$, $l$ and $r$ for
$\cat{Tricat}_2$. For the left unit constraint $l$, consider a lax homomorphism
$F \colon \S \to \T$, and write $F'$ for the composite $\id_\T.F \colon \S \to
\T$. Now, $F'$ agrees with $F$ on 0-cells and on hom-bicategories, but differs
in the remaining coherence data; indeed, we have
\begin{align*}
\iota^{F'}_{A} &= \cd{I_{FA} \ar@2[r]^-{\id_{I_{FA}}} & I_{FA} \ar@2[r]^-{\iota^F_{A}} &
FI_A}\\
\text{and} \quad \chi^{F'}_{f, g} &= \cd{Fg. Ff \ar@2[r]^-{\id_{Fg.Ff}} & Fg. Ff
\ar@2[r]^-{\chi^F_{f, g}} & F(gf).}
\end{align*}
Thus we define a 2-cell $l_F \colon \id_\T . F \Rightarrow F$ in
$\cat{Tricat}_2$ as follows:

\begin{enumerate}[(TD1)]
\item[(TD1)] $(l_F)_{A, B} = \id_{F_{A, B}} \colon {F_{A, B}} \Rightarrow
    {F_{A, B}}$;
\item[(TD2)] $M^{l_F}_A$ is the unit isomorphism $\iota^F_A.(\id_{I_{FA}})
    \Rrightarrow \iota^F_A$ in the bicategory $\T(FA, FA)$;
\item[(TD3)] $\Pi^{l_F}_{f, g}$ is the  unit isomorphism $\chi^F_{f,
    g}.(\id_{Fg.Ff}) \Rrightarrow \chi^F_{f, g}$ in the bicategory $\T(FA,
    FC)$.
\end{enumerate}

Now each axiom (TA1)--(TA4) is a tautology which describes how we obtained
$\chi^{F'}$, $\delta^{F'}$, $\gamma^{F'}$ and $\omega^{F'}$ from the
corresponding data for $F$. The definition of $r$ is dual to that of $l$, so we
pass over it and onto the associativity constraint $a$. Consider three lax
homomorphisms $F \colon \R \to \S$, $G \colon \S \to \T$ and $H \colon \T \to
\U$ and the two composites $(HG)F$ and $H(GF) \colon \R \to \U$. As above,
these agree on 0-cells and on hom-bicategories (and so we write their common
value simply as $HGF$) but differ with respect to coherence data. This time we
have:
\begin{align*}
\iota^{(HG)F} &= HG\iota^F.(H\iota^G.\iota^H)\text, & \iota^{H(GF)} &=
(HG\iota^F.H\iota^G).\iota^H\text,\\
\chi^{(HG)F} &= HG\chi^F.(H\chi^G.\chi^H)\qquad\text{and} & \chi^{H(GF)} &=
(HG\chi^F.H\chi^G).\chi^H\text,
\end{align*}
where we omit the subscripts for clarity. Thus we take $a_{F, G, H} \colon
(HG)F \Rightarrow H(GF)$ in $\cat{Tricat}_2$ to be:

\begin{enumerate}[(TD1)]
\item[(TD1)] $(a_{F,G,H})_{A, B} = \id_{(HGF)_{A, B}} \colon (HGF)_{A, B}
    \Rightarrow (HGF)_{A, B}$;
\item[(TD2)] $M^{a_{F,G,H}}_A$ is the associativity isomorphism
\[HG\iota^F_A.(H\iota^G_{FA}.\iota^H_{GFA}) \Rrightarrow
(HG\iota^F_A.H\iota^G_{FA}).\iota^H_{GFA}\] in the bicategory $\U(HGFA, HGFA)$;
\item[(TD3)] $\Pi^{a_{F,G,H}}_{f, g}$ is the  associativity isomorphism
\[HG\chi^F_{f,g}.(H\chi^G_{Ff,Fg}.\chi^H_{GFf,GFg}) \Rrightarrow
(HG\chi^F_{f,g}.H\chi^G_{Ff,Fg}).\chi^H_{GFf,GFg}\]
 in the bicategory $\U(HGFA, HGFC)$.
\end{enumerate}

We must now verify axioms (TA1)--(TA4) for these data. For this we observe that
the 3-cell data $\chi$,  $\gamma$, $\delta$ and $\omega$ for $H(GF)$ and for
$(HG)F$ are, in fact, obtained as different bracketings of the same pasting
diagram.
So by the
pasting theorem for bicategories, we can obtain the 3-cell data $\chi$,
$\gamma$, $\delta$ and $\omega$ for $H(GF)$ from that for $(HG)F$ by pasting
with suitable associativity isomorphisms in the appropriate hom-bicategory of
$\U$; and this is precisely what axioms (TA1)--(TA4) say.

It remains to check the naturality of $l$, $r$ and $a$, and the pentagon and
triangle identities. For the naturality of $l$, we must show that for any
2-cell $\alpha \colon F \Rightarrow G$ of $\cat{Tricat}_2$, the following
diagram commutes:
\[\cd{
    {\id_\T} . F \ar@2[r]^{l_F} \ar@2[d]_{\id_\T . \alpha} & F \ar@2[d]^\alpha \\
    {\id_\T} . G \ar@2[r]_{l_G} & G\text.}
\]
We easily verify that the left-hand 2-cell $\alpha' = \id_\T . \alpha$ has
components $\alpha'_\theta = \alpha_\theta$, $M^{\alpha'}_A =
M^{\alpha}_A.(\id_{I_{FA}})$ and $\Pi^{\alpha'}_{f, g} = \Pi^{\alpha}_{f, g} .
(\id_{Ff. Fg})$; therefore the naturality of $l$ is a consequence of the
naturality of the left unit constraints in the hom-bicategories of $\T$; and
dually for $r$. For the naturality of $a$, we must show that the following
diagram commutes in $\cat{Tricat}_2$ for all suitable 2-cells $\alpha$, $\beta$
and $\epsilon$:
\[\cd{
    (HG)F \ar@2[r]^{a_{F,G,H}} \ar@2[d]_{(\alpha\beta)\epsilon} & H(GF) \ar@2[d]^{\alpha(\beta\epsilon)} \\
    (H'G')F' \ar@2[r]_{a_{F',G',H'}} & H'(G'F')\text,}
\]
for which we must show that (TD1)--(TD3) agree for the two ways around this
square. For (TD1) this is trivial; so consider (TD2). For both $(\alpha
\beta)\epsilon$ and $\alpha(\beta\epsilon)$, we obtain this datum by pasting
together the same $3 \times 3$ diagram of 3-cells; the only difference being
the manner in which we bracket together the boundary of this diagram. Thus the
commutativity of the above square with respect to (TD2) is a further instance
of the pasting theorem for bicategories. (TD3) is obtained in a similar manner.

Finally, it is not hard to verify that the pentagon and triangle identities for
$a$, $l$ and $r$ follow from instances of the pentagon and triangle identities
in the hom-bicategories of the target tricategory. This completes the
definition of the bicategory $\cat{Tricat}_2$.

\newcommand{\smt}{\mathbin{{}_s \times_t}}
\newcommand{\DI}[1]{I_{#1}}

\section{Towards a tricategory of tricategories}\label{Sec:towards}
We now wish to describe a tricategory of tricategories $\cat{Tricat}_3$. This
will have the same 0-cells and 1-cells as $\cat{Tricat}_2$, but will have
2-cells with one fewer level of degeneracy, which consequently admit a notion
of 3-cell between them. Although we introduce the 2-~and 3-cells of
$\cat{Tricat}_3$ in this Section, we will not actually prove that we obtain a
tricategory from them until we reach Section~\ref{tritri}. As explained in the
Introduction, we do this for two reasons. Firstly, so that we can set up some
machinery which will allow us to avoid checking all the tricategorical
coherence axioms by hand; and secondly, in order to investigate the curious
fact that $\cat{Tricat}_3$ does not really extend $\cat{Tricat}_2$, in that not
every $2$-cell of the latter gives rise to a $2$-cell of the former.

We now begin our description of $\cat{Tricat}_3$. Its objects and $1$-cells
are, as stated above, tricategories and lax trihomomorphisms. The $2$-cells are
to be ``singly degenerate oplax tritransformations''. The most obvious way of
interpreting this notion would be as follows:
\begin{Defn}
Let there be given lax homomorphisms of tricategories $F, G \colon \S \to \T$;
then an
\defn{oplax icon} $\alpha \colon F \Tor G$ may exist only if $F$ and $G$ agree on objects
whereupon it consists of the following data:
\begin{enumerate}[(ID1)]
\item[(ID1)] For each $A$ and $B$ in $\S$, an oplax natural transformation
\[\alpha_{A, B} \colon F_{A, B} \Rightarrow G_{A, B} \colon \S(A, B) \to \T(FA,
FB)\] (so in particular, for each 1-cell $f \colon A \to B$ of $\S$, we have a 2-cell $\alpha_f
\colon Ff \Rightarrow Gf$ of $\T$, and for each 2-cell $\theta \colon f
\Rightarrow g$ of $\S$, a 3-cell
\[\cd{
    Ff \ar@2[r]^{F \theta} \ar@2[d]_{\alpha_f} \dthreecell{dr}{\alpha_{\theta}}& Fg \ar@2[d]^{\alpha_g} \\
    Gf \ar@2[r]_{G \theta} & Gg\quad \text{);}
}\]
\item[(ID2)] For each object $A$ of $\S$, a 3-cell of $\T$:
\[
\cd{
    I_{FA} \ar@2[r]^{\iota^F_A} \ar@{=}[d] \dthreecell{dr}{M^{\alpha}_A} & FI_A \ar@2[d]^{\alpha_{I_A}} \\
    I_{GA} \ar@2[r]_{\iota^G_A} & GI_A\text;
}\]
\item[(ID3)] For each $A, B$ and $C$ in $\S$, a modification
\[
\cd{
  F(\thg)\otimes F(?) \ar@2[r]^{\chi^F} \ar@2[d]_{\alpha_{(\thg)}\otimes \alpha_{(?)}} \dthreecell{dr}{\Pi_{A, B, C}^\alpha} &
  F\big((\thg) \otimes (?)\big) \ar@2[d]^{\alpha_{((\thg) \otimes (?))}} \\
  G(\thg)\otimes G(?) \ar@2[r]_{\chi^G} &
  G\big((\thg) \otimes (?)\big)\text,
}
\]
where, for instance, $F(\thg) \otimes F(?)$ represents the homomorphism
\[\S(B, C) \times \S(A, B) \xrightarrow{F \times F} \T(FB, FC) \times \T(FA, FB) \xrightarrow{\otimes} \T(FA,
FC)\] (so in particular, for each pair of composable 1-cells $f \colon A \to B$, $g \colon B \to
C$ of $\S$, we have a 3-cell of $\T$:
\[
\cd{
    Fg . Ff \ar@2[r]^{\chi^F_{f, g}} \ar@2[d]_{\alpha_g . \alpha_f} \dthreecell{dr}{\Pi^{\alpha}_{f, g}} & F(g f) \ar@2[d]^{\alpha_{gf}} \\
    Gg . Gf \ar@2[r]_{\chi^G_{f, g}} & G(g f)
}\quad\text{).}\]
\end{enumerate}
These data are subject to the following axioms:
\begin{enumerate}[(IA1)]
\item[(IA1)] For each 1-cell $f \colon A \to B$ of $\S$, the following
    pasting equality holds:
\[
 \hphantom{=} \quad \cd[@!C@!R@R+1em]{
    & FI_{B}.Ff
      \ar@2[r]^{\chi^F} \dthreecell{ddr}{\gamma^F}
    & F(I_B.f)
      \ar@2[dr]^{F\l}
    \\ I_{FB}.Ff
      \ar@2[dr]^{\l} \ar@2[ur]^{\iota^F.1} \ar@2[d]_{1 . \alpha_f}
    & & & Ff
      \ar@2[d]^{\alpha_f}
    \\ I_{GB}.Gf
      \ar@2[dr]_{\l} \twocong{r}
    & Ff
      \ar@{=}[r] \ar@2[d]^{\alpha_f} \twoeq{dr}
    & Ff
      \ar@{=}[ur] \ar@2[d]^{\alpha_f} \twoeq{r}
    & Gf
    \\ & Gf
      \ar@{=}[r]
    & Gf
      \ar@{=}[ur]
}\]
\[
= \quad \cd[@!C@!R@R+1em]{
    & FI_B . Ff
      \ar@2[r]^{\chi^F} \ar@2[d]|{\alpha_{I_B} . \alpha_f} \dthreecell{dr}{\Pi^\alpha}
    & F(I_B . f)
      \ar@2[dr]^{F\l} \ar@2[d]|{\alpha_{I_B . f}}
    \\ I_{FB}.Ff
      \ar@2[ur]^{\iota^F.1} \ar@2[d]_{1 . \alpha_f} \ar@{}[r] \ar@3?(0.5)+/u  0.2cm/;?(0.5)+/d 0.2cm/_(1.5)*+[l]{\scriptstyle \ \ \overline{M^\alpha.1}}
    & GI_B . Gf
      \ar@2[r]_{\chi^G} \dthreecell{ddr}{\gamma^G}
    & G(I_B . f)
      \ar@2[dr]_{G \l} \dthreecell{r}{\alpha_\l}
    & Ff
      \ar@2[d]^{\alpha_f}
    \\ I_{GB}.Gf
      \ar@2[dr]_{\l} \ar@2[ur]_{\iota^G.1}
    & & & Gf\text.
    \\ & Gf
      \ar@{=}[r]
    & Gf
      \ar@{=}[ur]
}\]
\item[(IA2)] For each 1-cell $f \colon A \to B$ of $\S$, the following
    pasting equality holds:
\[
 \hphantom{=} \quad \cd[@!C@!R@R+1em]{
    & Ff.FI_{A}
      \ar@2[r]^{\chi^F} \dthreecell{ddr}{\delta^F}
    & F(f.I_A)
      \ar@2[dr]^{F\r}
    \\ Ff.I_{FA}
      \ar@2[dr]^{\r} \ar@2[ur]^{1.\iota^F} \ar@2[d]_{\alpha_f . 1}
    & & & Ff
      \ar@2[d]^{\alpha_f}
    \\ Gf.I_{GA}
      \ar@2[dr]_{\r} \twocong{r}
    & Ff
      \ar@{=}[r] \ar@2[d]^{\alpha_f} \twoeq{dr}
    & Ff
      \ar@{=}[ur] \ar@2[d]^{\alpha_f} \twoeq{r}
    & Gf
    \\ & Gf
      \ar@{=}[r]
    & Gf
      \ar@{=}[ur]
}\]
\[
= \quad \cd[@!C@!R@R+1em]{
    & Ff.FI_A
      \ar@2[r]^{\chi^F} \ar@2[d]|{\alpha_f . \alpha_{I_A}} \dthreecell{dr}{\Pi^\alpha}
    & F(f.I_A)
      \ar@2[dr]^{F\r} \ar@2[d]|{\alpha_{f . I_A}}
    \\ Ff.I_{FA}
      \ar@2[ur]^{1.\iota^F} \ar@2[d]_{\alpha_f . 1} \ar@{}[r] \ar@3?(0.5)+/u  0.2cm/;?(0.5)+/d 0.2cm/_(1.5)*+[l]{\scriptstyle \ \ \overline{1.M^\alpha}}
    & Gf.GI_A
      \ar@2[r]_{\chi^G} \dthreecell{ddr}{\delta^G}
    & G(f.I_A)
      \ar@2[dr]_{G \r} \dthreecell{r}{\alpha_\r}
    & Ff
      \ar@2[d]^{\alpha_f}
    \\ Gf.I_{GA}
      \ar@2[dr]_{\r} \ar@2[ur]_{1.\iota^G}
    & & & Gf\text.
    \\ & Gf
      \ar@{=}[r]
    & Gf
      \ar@{=}[ur]
}\]
\item[(IA3)] For each triple $f, g, h$ of composable 1-cells of $\S$, the
    following pasting equality holds:
\[
 \hphantom{=} \quad \cd[@!C@!R@R+1em]{
    & F(hg). Ff
      \ar@2[r]^{\chi^F} \dthreecell{ddr}{\omega^F}
    & F((hg)f)
      \ar@2[dr]^{F\a}
    \\ (Fh.Fg).Ff
      \ar@2[dr]^{\a} \ar@2[ur]^{\chi^F.1} \ar@2[d]_{(\alpha_h .
      \alpha_g) . \alpha_f}
    & & & F(h(gf))
      \ar@2[d]^{\alpha_{h(gf)}}
    \\ (Gh.Gg).Gf
      \ar@2[dr]_{\a} \twocong{r}
    & Fh.(Fg. Ff)
      \ar@2[r]^{1. \chi^F} \ar@2[d]|{\alpha_h . (\alpha_g . \alpha_f)} \dthreecell{dr}{\overline
      {1.\Pi^\alpha}}
    & Fh. F(gf)
      \ar@2[ur]^{\chi^F} \ar@2[d]^{\alpha_h . \alpha_{gf}} \dthreecell{r}{\Pi^\alpha}
    & G(h(gf))
    \\ & Gh.(Gg.Gf)
      \ar@2[r]_{1.\chi^G}
    & Gh.G(gf)
      \ar@2[ur]_{\chi^G}
}\]
\[
= \quad \cd[@!C@!R@R+1em]{
    & F(hg). Ff
      \ar@2[r]^{\chi^F} \ar@2[d]^{\alpha_{hg} . \alpha_f} \dthreecell{dr}{\Pi^\alpha}
    & F((hg)f)
      \ar@2[dr]^{F\a} \ar@2[d]^{\alpha_{(hg)f}}
    \\ (Fh.Fg).Ff
      \ar@2[ur]^{\chi^F.1} \ar@2[d]_{(\alpha_h . \alpha_g) . \alpha_f} \ar@{}[r] \ar@3?(0.5)+/u  0.2cm/;?(0.5)+/d 0.2cm/_(1.5)*+[l]{\scriptstyle \ \ \overline{\Pi^\alpha.1}}
    & G(hg).Gf
      \ar@2[r]_{\chi^G} \dthreecell{ddr}{\omega^G}
    & G((hg)f)
      \ar@2[dr]_{G \a} \dthreecell{r}{\alpha_\a}
    & F(h(gf))
      \ar@2[d]^{\alpha_{h(gf)}}
    \\ (Gh.Gg).Gf
      \ar@2[dr]_{\a} \ar@2[ur]_{\chi^G.1}
    & & & G(h(gf))\text.
    \\ & Gh.(Gg.Gf)
      \ar@2[r]_{1.\chi^G}
    & Gh.G(gf)
      \ar@2[ur]_{\chi^G}
}\]
\end{enumerate}
\end{Defn}
The definition of oplax icon generalises that of ico-icon, in that:
\begin{Prop}\label{embed} Let $F, G \colon \S \to \T$ be lax
homomorphisms. Then the ico-icons \mbox{$F \Rightarrow G$} are in bijection
with the class of oplax icons $\alpha \colon F \Tor G$ for which each component
$\alpha_f \colon Ff \Rightarrow Gf$ is an identity 2-cell: they are the
\emph{\underline{i}}dentity \emph{\underline{c}}omponents
\emph{\underline{o}}plax icons.
\end{Prop}

Unfortunately, oplax icons do not provide a suitable notion of 2-cell for our
tricategory $\cat{Tricat}_3$. The reason is that although oplax icons may be
``whiskered'' with lax homomorphisms on each side, these whiskerings do not
give rise to a well-defined composition of oplax icons along a 0-cell
boundary. Indeed, if we are given a diagram
\[\cd{
  {\S} \ar@/^1em/[r]^{F} \ar@/_1em/[r]_{F'} \dtwocell{r}{\alpha} &
  {\T} \ar@/^1em/[r]^{G} \ar@/_1em/[r]_{G'} \dtwocell{r}{\beta} &
  {\U}
}\] of lax homomorphisms and oplax icons, then there are two canonical ways of
composing it up which need not agree, even up to isomorphism. The same
phenomenon occurs if one tries to form a tricategory of bicategories whose
2-cells are oplax natural transformations. In order to obtain a tricategory, we
therefore restrict attention to a suitable subclass of the oplax icons:
\begin{Defn}\label{pseudoicon}
Let $F, G \colon \S \to \T$ be lax homomorphisms. By a \defn{pseudo-icon}
\mbox{$\alpha \colon F \Tor G$} we mean an oplax icon
    $\alpha$ for which each 3-cell $\alpha_\theta$, $M^\alpha_A$, and
    $\Pi^\alpha_{f, g}$ is invertible.
\end{Defn}

These pseudo-icons are to be the 2-cells of $\cat{Tricat}_3$. Note that,
although every ico-icon gives rise to an oplax icon, it is only the
\emph{invertible} ico-icons which give rise to pseudo-icons. We now turn to the
3-cells of $\cat{Tricat}_3$.
\begin{Defn}\label{pseudoiconmod}
Given pseudo-icons $\alpha$, $\beta \colon F \Tor G$, a
\defn{pseudo-icon modification} \mbox{$\Gamma \colon \alpha \Rrightarrow
\beta$} consists in the following data:
\begin{enumerate}[(MD1)]
\item[(MD1)] For each $A, B$ in $\S$, a modification $\Gamma_{A, B} \colon
    \alpha_{A, B} \Rrightarrow \beta_{A, B}$ (and so in particular, for
    each
 1-cell $f \colon A \to B$ of $\S$, a
3-cell $\Gamma_f \colon \alpha_f \Rrightarrow \beta_f$ of $\T$);
\end{enumerate}
subject to the following axioms:
\begin{enumerate}[(MA1)]
\item[(MA1)] For each object $A$ of $\S$, the following pasting equality
    holds:
\[\cd[@u@!@+2em]{
    I_{FA} \ar@{=}[r] \ar@{=}[d] \ar@2[dr]_{\iota^{F}_A} &
    I_{FA} \ar@2[dr]^{\iota^F_A}  \\
    I_{GA} \ar@2[dr]_{\iota^{G}_A} \dlthreecell{r}{M^{\beta}_A}  &
    FI_A \ar@2[d]^{\beta_{I_A}} \ar@{=}[r] \dthreecell{dr}{\Gamma_{I_A}} &
    FI_A \ar@2[d]^{\alpha_{I_A}}\\
    & GI_A \ar@{=}[r] & GI_A
}\quad = \quad \cd[@u@!@+2em]{
    I_{FA} \ar@{=}[r] \ar@{=}[d] &
    I_{FA} \ar@{=}[d] \ar@2[dr]^{\iota^F_A} \\
    I_{GA} \ar@{=}[r] \ar@2[dr]_{\iota^{G}_A}  &
    I_{GA} \ar@2[dr]_{\iota^G_A} \dlthreecell{r}{M^{\alpha}_A}  &
    FI_A \ar@2[d]^{\alpha_{I_A}}\\
    & GI_A \ar@{=}[r] & GI_A\text;
}
\]
\item[(MA2)] For each pair of composable 1-cells $f \colon A \to B$, $g
    \colon B \to C$ of $\S$, the following pasting equality holds:
\end{enumerate}
\[\cd[@u@!@+2em@R-1em]{
    Fg . Ff \ar@{=}[r] \ar@2[d]_{\beta_g . \beta_f} \ar@2[dr]_{\chi^{F}_{f, g}} &
    Fg . Ff \ar@2[dr]^{\chi^F_{f, g}}  \\
    Gg . Gf \ar@2[dr]_{\chi^{G}_{f, g}} \dlthreecell{r}{\Pi^{\beta}_{f, g}}  &
    F(gf) \ar@2[d]^{\beta_{gf}} \ar@{=}[r] \dthreecell{dr}{\Gamma_{gf}} &
    F(gf) \ar@2[d]^{\alpha_{gf}}\\
    & G(gf) \ar@{=}[r] & G(gf)
}\quad = \quad \cd[@u@!@+2em@R-1em]{
    Fg . Ff \ar@{=}[r] \ar@2[d]_{\beta_g . \beta_f} \dthreecell{dr}{\Gamma_g . \Gamma_f} &
    Fg . Ff \ar@2[d]^{\alpha_g . \alpha_f} \ar@2[dr]^{\chi^F_{f, g}} \\
    Gg . Gf \ar@{=}[r] \ar@2[dr]_{\chi^{G}_{f, g}}  &
    Gg . Gf \ar@2[dr]_{\chi^G_{f, g}} \dlthreecell{r}{\Pi^{\alpha}_{f, g}}  &
    F(gf) \ar@2[d]^{\alpha_{gf}}\\
    & G(gf) \ar@{=}[r] & G(gf)\text.
}
\]
\end{Defn}
\begin{Thm}\label{mainresultprop}
There is a tricategory $\cat{Tricat}_3$ with objects being tricategories;
\mbox{1-cells}, lax homomorphisms; 2-cells, pseudo-icons; and 3-cells,
pseudo-icon modifications.
\end{Thm}
It would certainly be possible to prove this result at this point in the paper:
we would simply follow the same path as in Section~\ref{bicat-tricat}, first
defining the various  kinds of composition we need, then the various pieces of
coherence data, and finally checking the coherence axioms these must satisfy.
However, rather than doing this directly, we would like to reuse some of the
results we proved about $\cat{Tricat}_2$.

Indeed, we have already shown that that the composition of lax homomorphisms is
associative up to an invertible ico-icon. Each invertible ico-icon witnessing this
associativity gives rise to a corresponding pseudo-icon in $\cat{Tricat}_3$;
and so by taking these pseudo-icons as our witnesses for associativity in
$\cat{Tricat}_3$, we might hope to be able to reuse the coherence work we did
in Section~\ref{bicat-tricat}.

However, matters are not quite this simple. If we take the unique sensible
definition of vertical composition of pseudo-icons, then we find that the
composition of two invertible ico-icons \emph{qua} ico-icon does not agree with
their composite \emph{qua} pseudo-icon. In particular, the invertible ico-icons
witnessing associativity in $\cat{Tricat}_2$ become mere equivalence
pseudo-icons in $\cat{Tricat}_3$, whilst each commutative diagram of coherence
2-cells in $\cat{Tricat}_2$ gives rise to a diagram in $\cat{Tricat}_3$ which
may commute only up to an invertible 3-cell.

Intuitively, it is clear that this should not be a problem, and that we should
still be able to ``read off'' the coherence for $\cat{Tricat}_3$ from that for
$\cat{Tricat}_2$, but to make this precise we must turn our intuition into a
mathematical principle. In order to motivate how we will do this, let us
examine more closely why the naive approach does not work.

The problem is essentially that the putative tricategory $\cat{Tricat}_3$ does not include
\emph{all} of the data carried by the mere bicategory $\cat{Tricat}_2$. This occurs at the
level of basic cell data~--~since not every ico-icon is a pseudo-icon~--~but more importantly,
at the level of compositional data: the data for the strictly associative composition of
ico-icons from $\cat{Tricat}_2$ is no longer present in $\cat{Tricat}_3$.

The solution we give to this problem is to describe a categorical structure into which
tricategories, lax homomorphisms, pseudo-icons and modifications may be formed which is richer
than $\cat{Tricat}_3$, and in particular retains \emph{all} the data from $\cat{Tricat}_2$.
This categorical structure is not a tricategory, but rather what we call a \emph{locally
cubical bicategory}. This is a genuinely weak three-dimensional structure whose coherence laws
are particularly simple: they have a bicategorical rather than tricategorical flavour. In
particular, the locally cubical bicategory of tricategories that we construct will be able to
take its coherence data directly from $\cat{Tricat}_2$.

The existence of the desired tricategory of tricategories $\cat{Tricat}_3$ now
follows from a general result (given in Section~\ref{Sec:fromltt}) which says
that any well-behaved locally cubical bicategory gives rise to a tricategory in
a canonical way. This result can be seen as a crystallisation of the intuition
we had above that we should be able to ``read off'' the coherence of
$\cat{Tricat}_3$ from $\cat{Tricat}_2$.

%
%
%
%
%
%
%
%
%
%
%
%

\section{Locally cubical bicategories}\label{Sec:locallydouble}
The purpose of this section is to define the locally cubical bicategories
alluded to at the end of the previous section. Like tricategories, these are
weak categorical structures comprised of 0-, 1-, 2-~and 3-cells; however, the
2-cells come in two varieties, \emph{horizontal} and \emph{vertical}, whilst
the 3-cells are cubical in nature. Composition of vertical 2-cells is strictly
associative; that of horizontal 2-cells is only so up to an invertible 3-cell;
whilst the associativity constraints for 1-cells are given by \emph{vertical}
2-cells, and are of a bicategorical, up-to-isomorphism, rather than a
tricategorical, up-to-equivalence, kind. A locally cubical bicategory may be
described succinctly as a ``bicategory weakly enriched in pseudo double
categories''; and our task in this section will be to expand upon this
description.

 The concept of \emph{strict} double
category is due to Ehresmann. It is an example of the notion of \emph{double
model} for an essentially-algebraic theory, this being a model of the theory in
its own category of ($\cat{Set}$-based) models. Thus a double category~--~which
is a double model of the theory of categories~--~is a category object in
$\cat{Cat}$.

 The theory
of categories is somewhat special, since its category of ($\cat{Set}$-based)
models may be enriched to a 2-category, so that, as well as \emph{strict}
category objects in $\cat{Cat}$, we may also consider \emph{pseudo} category
objects: and these are the pseudo double categories which we will be interested
in.

\begin{Defn}
A \defn{pseudo double category}~\cite{GP1} $\mathfrak C$ is given by specifying
a collection of \emph{objects} $x, y, z, \dots$, a collection of \emph{vertical
1-cells} between objects, which we write as $a \colon x \to y$, a collection of
\emph{horizontal 1-cells} between objects, which we write as $f \colon x \tor
y$, and a collection of \emph{2-cells}, each of which is bounded by a square of
horizontal and vertical arrows, and which we write as:
\[\cd{
x \ar[d]_{a} \ar[r]|-{\object@{|}}^{f} \dtwocell{dr}{\alpha}& w \ar[d]^{b}  \\
y \ar[r]|-{\object@{|}}_{g} & z\text,}\] or sometimes simply as $\alpha \colon f \Rightarrow g$.
Moreover, we must give:
\begin{itemize}
\item Identities and composition for vertical 1-cells, $\id_x \colon x \to
    x$ and $(a, b) \mapsto ab$, making the objects and vertical arrows into
    a category $\C_0$;
\item Vertical identities and composition for 2-cells, $\id_{f} \colon f
    \Rightarrow f$ and $(\beta, \alpha) \mapsto \beta \alpha$:
\[\cd{
x \ar[d]_{\id_x} \ar[r]|-{\object@{|}}^{f} \dtwocell{dr}{\id_{f}}& y \ar[d]^{\id_y}  \\
x \ar[r]|-{\object@{|}}_{f} & y} \quad \text; \quad \cd{
u \ar[d]_{a} \ar[r]|-{\object@{|}}^{f} \dtwocell{dr}{\alpha}& x \ar[d]^{b}  \\
v \ar[r]|-{\object@{|}}_{g} \ar[d]_{c} \dtwocell{dr}{\beta}& y  \ar[d]^{d} \\
w \ar[r]|-{\object@{|}}_{h} & z} \mapsto
\cd{u \ar[d]_{ca} \ar[r]|-{\object@{|}}^{f} \dtwocell{dr}{\beta\alpha}& x \ar[d]^{db}  \\
w \ar[r]|-{\object@{|}}_{h} & z}
\]
making the horizontal arrows and 2-cells into a category $\C_1$ for which
``vertical source'' and ``vertical target'' become functors $s, t \colon
\C_1 \to \C_0$;
\item Identities and composition for horizontal 1-cells, $\DI x \colon x
    \tor x$ and $(g, f) \mapsto gf$;
\item Horizontal identities and composition for 2-cells, $\DI a \colon \DI
    x \Rightarrow \DI y$ and $(\beta, \alpha) \mapsto \beta \ast \alpha$:
\[\cd{
x \ar[d]_{a} \ar[r]|-{\object@{|}}^{\DI x} \dtwocell{dr}{\DI a}& x \ar[d]^{a}  \\
y \ar[r]|-{\object@{|}}_{\DI y} & y} \quad \text; \quad  \cd{
u \ar[d]_{a} \ar[r]|-{\object@{|}}^{f} \dtwocell{dr}{\alpha}& v \ar[d]^{b} \ar[r]|-{\object@{|}}^{g} \dtwocell{dr}{\beta} & w \ar[d]^{c} \\
x \ar[r]|-{\object@{|}}_{h} & y \ar[r]|-{\object@{|}}_{k} & z} \mapsto
\cd{u \ar[d]_{a} \ar[r]|-{\object@{|}}^{gf} \dtwocell{dr}{\beta \ast \alpha}& w \ar[d]^{c}  \\
x \ar[r]|-{\object@{|}}_{kh} & z\text,}
\]
satisfying functoriality constraints: firstly, $\DI{(\thg)}$ is a functor
$\C_0 \to \C_1$, which says that $\DI{\id_x} = \id_{\DI x}$ and $\DI{ab} =
\DI a.\DI b$ and secondly, horizontal composition is a functor
$\mathord{\ast} \colon \C_1 \smt \C_1 \to \C_1$ which says that $\id_{g}
\ast \id_{f} = \id_{gf}$ and $(\delta \ast \gamma).(\beta \ast \alpha) =
(\delta \beta) \ast (\gamma \alpha)$.
\item Horizontal unitality and associativity constraints given by 2-cells
\[\cd{
x \ar[d]_{\id_x} \ar[r]|-{\object@{|}}^{\DI y.f} \dtwocell{dr}{\l_{f}}& y \ar[d]^{\id_y}  \\
x \ar[r]|-{\object@{|}}_{f} & y}\text, \quad \cd{
x \ar[d]_{\id_x} \ar[r]|-{\object@{|}}^{{f}.\DI x} \dtwocell{dr}{\r_{f}}& y \ar[d]^{\id_y}  \\
x \ar[r]|-{\object@{|}}_{f} & y}\text, \quad \text{and} \quad \cd[@C+1.5em]{
x \ar[d]_{\id_x} \ar[r]|-{\object@{|}}^{h(gf)} \dtwocell{dr}{\a_{f,g,h}}& z \ar[d]^{\id_z}  \\
x \ar[r]|-{\object@{|}}_{(hg)f} & z\text,}
\]
natural in $f$, $g$ and $h$, and invertible as arrows of $\C_1$. These
2-cells must obey two laws: the pentagon law, which equates the two routes
from $k(h(gf))$ to $((kh)g)f$, and the triangle law, which equates the two
routes from $g.(\DI y.f)$ to $gf$.
\end{itemize}
\end{Defn}
Pseudo double categories are sometimes also known as \emph{weak} double
categories; they are a special case of Verity's more general notion of
\emph{double bicategory}~\cite{dv}. A more comprehensive reference on pseudo
double categories is \cite{GP1}: though be aware that we interchange its usage
of the terms ``horizontal'' and ``vertical'' to give a better fit with the
usual \mbox{2-categorical} terminology. Since the only sorts of double
categories we will be concerned with in this paper are the pseudo ones, we may
sometimes choose to write simply ``double category'', leaving the qualifier
``pseudo'' understood.

Some simple examples of pseudo double categories are $\mathfrak{Cat}$, the pseudo double
category of ``categories, functors, profunctors and transformations'', $\mathfrak{Rng}$, the
pseudo double category of ``rings, ring homomorphisms, bimodules and skew-linear maps'', and
the pseudo double category $\mathfrak{Span}(\C)$ of ``objects, morphisms, spans and span
morphisms'' in a category with pullbacks $\C$. These are typical examples of pseudo double
categories, in that they have notions of \emph{homomorphism} and \emph{bimodule} as their
respective vertical and horizontal 1-cells. Any bicategory $\B$ gives us a pseudo double
category $\UU(\B)$ with only identity vertical 1-cells, whilst any pseudo double category
$\mathfrak C$ gives a bicategory $\HH(\mathfrak C)$ upon throwing away the non-identity
vertical 1-cells, and all the 2-cells except for those whose vertical source and target are
identity arrows. We will refer to such 2-cells as \defn{globular 2-cells}; they are also
sometimes known as \emph{special} 2-cells.

Just as in the theory of bicategories, the appropriate notion of morphism
between pseudo double categories only preserves horizontal composition up to
comparison 2-cells, the most important case being the \emph{homomorphisms}, for
which these 2-cells are invertible. We can define a homomorphism between small
pseudo double categories in terms of a pseudomorphism of pseudocategory
objects, but just as easy is to give the elementary description:
\begin{Defn}
A \defn{homomorphism of pseudo double categories} $F \colon \mathfrak C \to
\mathfrak D$ is given by assignations on objects, 1-cells and 2-cells which
preserve source and target and are functorial with respect to vertical
composition of 1-\ and 2-cells, together with comparison 2-cells
\[\cd{
{Fx} \ar[d]_{\id_{F x}} \ar[r]|-{\object@{|}}^{\DI{Fx}} \dtwocell{dr}{\mathfrak m_{x}}& {F x} \ar[d]^{\id_{F x}}  \\
Fx \ar[r]|-{\object@{|}}_{F\DI x} & Fx} \quad \text{and} \quad \cd{
{F x} \ar[d]_{\id_{F x}} \ar[r]|-{\object@{|}}^{F g. F f} \dtwocell{dr}{\mathfrak m_{f, g}}& {F z} \ar[d]^{\id_{F z}}  \\
Fx \ar[r]|-{\object@{|}}_{F (g f)} & Fz}
\]
which are invertible as arrows of $\D_1$, and natural in $x$, respectively $g$
and $f$. Moreover, we require the commutativity of three familiar diagrams,
which equate, respectively, the two possible ways of going from $Ff.\DI{Fx}$ to
$Ff$, from $\DI{Fy}. Ff$ to $Ff$, and from $Fh.(Fg.Ff)$ to $F((hg)f)$.
\end{Defn}
With the obvious notion of composition and identities, we obtain a category
$\cat{DblCat}$ of (possibly large) pseudo double categories and homomorphisms
between them. If we write $\cat{Bicat}$ for the category of bicategories and
homomorphisms, then the assignations $\B \mapsto \UU(\B)$ and $\mathfrak C
\mapsto \HH(\mathfrak C)$ described above extend to a pair of adjoint functors
\mbox{$\UU \dashv \HH \colon \cat{DblCat} \to \cat{Bicat}$}, for which the
composite $\HH\UU$ is the identity; we can thus view $\cat{Bicat}$ as a
coreflective subcategory of $\cat{DblCat}$.

Now, $\cat{DblCat}$ is in fact the underlying ordinary category of a 2-category
whose 2-cells are the so-called \emph{vertical transformations}. We can
understand these 2-cells by observing that there is a 2-monad  on the
2-category \mbox{$\cat{CatGph} := [\,\bullet \rightrightarrows \bullet, \,
\cat{Cat}]$} whose strict algebras are small pseudo double categories, and
whose algebra pseudomorphisms are the homomorphisms between them. The
corresponding \emph{algebra 2-cells} are precisely the vertical
transformations. Spelling this out, we have:
\begin{Defn}
A \defn{vertical transformation} $\alpha \colon F \Rightarrow G$ between
homomorphisms of pseudo double categories $F, G \colon \mathfrak C \to
\mathfrak D$ is given by specifying, for each object $x \in \mathfrak C$, a
vertical 1-cell $\alpha_x \colon Fx \to Gx$ of $\mathfrak D$ and for each
horizontal 1-cell $f \colon x \tor y$ in $\mathfrak C$ a 2-cell
\[\cd{
Fx \ar[d]_{\alpha_x} \ar[r]|-{\object@{|}}^{Ff} \dtwocell{dr}{\alpha_f}& Fy \ar[d]^{\alpha_y}  \\
Gx \ar[r]|-{\object@{|}}_{Gf} & Gy}\] of $\mathfrak D$, such that the $\alpha_x$'s are natural in
morphisms of $\D_0$, the $\alpha_f$'s are natural in morphisms of $\D_1$, and
the following diagrams commute:
\[\cd{
\DI{Fx} \ar@2[r]^{\mathfrak m^F_x} \ar@2[d]_{\DI{\alpha_x}} & F\DI x \ar@2[d]^{\alpha_{\DI x}} \\
\DI{Gx} \ar@2[r]_-{\mathfrak m^G_x} & G\DI x
 } \quad \text{and} \quad
\cd{ Fg.Ff \ar@2[r]^{\mathfrak m^F_{g, f}} \ar@2[d]_{\alpha_g \ast \alpha_f} & F(gf) \ar@2[d]^{\alpha_{gf}} \\
Gg.Gf \ar@2[r]_-{\mathfrak m^G_{g, f}} & G(gf)\text.}
 \]
\end{Defn}
In the case that $\mathfrak C$ and $\mathfrak D$ are bicategories, the vertical
transformation between homomorphisms $\mathfrak C \to \mathfrak D$ are
precisely the bicategorical icons of Section 1; however, the reader should
carefully note that the coreflection of $\cat{DblCat}$ into $\cat{Bicat}$ does
\emph{not} enrich to a two-dimensional coreflection, since there is no way of
coreflecting a general vertical transformation between homomorphisms of pseudo
double categories into an icon between the corresponding homomorphisms of
bicategories.

It follows from the algebraic description of $\cat{DblCat}$ that it admits a
wide class of 2-dimensional limits, of which we will only be concerned with
finite products. That $\cat{DblCat}$ admits these, makes it, of course, into a
symmetric monoidal category, but the 2-dimensional aspect of these products
means that we may view it instead as a \emph{symmetric monoidal 2-category}:
that is, a symmetric monoidal category whose tensor product is a 2-functor and
whose coherence natural transformations are 2-natural transformations. What we
now wish to describe is how we can use this monoidal 2-category $\cat{DblCat}$
as a suitable base for \emph{enrichment}.

For any monoidal category $\V$, we have the well-known notion of a \emph{category enriched in
$\V$} or \emph{$\V$-category}, which instead of having hom-sets between 0-cells, has
hom-objects drawn from $\V$, with the corresponding composition being expressed by morphisms
of $\V$ subject to associativity and unitality laws. Now, if instead of a monoidal category
$\V$ we begin with a \emph{monoidal bicategory} $\W$ in the sense of \cite{monoidalbicats},
then we may generalise this definition to obtain the notion of \emph{bicategory enriched in
$\W$} or \mbox{\emph{$\W$-bicategory}}~\cite{Carmody, Lack}. A $\W$-bicategory is like a
bicategory, but instead of hom-categories between 0-cells, it has hom-objects drawn from $\W$:
and instead of composition functors, it has composition morphisms drawn from $\W$, which are
now required to be associative and unital only up to coherent 2-cells of $\W$.\footnote{Note
that this differs from the notion of ``category enriched in a bicategory'' studied in
\cite{Enrichedbicats}; these are the \emph{polyads} of \cite{bicats}, and are essentially
categories enriched in a monoidal category where that monoidal category happens to be spread
out over many objects.}\ \  Thus we can think of a $\W$-bicategory as being a ``category
weakly enriched in $\W$''.

The simplest sort of enriched bicategory is a $\cat{Cat}$-bicategory, which is
just a (\mbox{locally} small) ordinary bicategory. Other examples are obtained
by taking \mbox{$\W = \cat{\V}\text-\cat{Cat}$} for some monoidal category
$\V$, for which a $\W$-bicategory has sets of 0-\ and 1-\ cells as usual, but
now a $\V$-object of 2-cells between any parallel pair of 1-cells; by taking
$\W = \cat{Mod}$, the bicategory of categories and profunctors, for which a
$\W$-bicategory is a \emph{probicategory} in the sense of Day \cite{Biclosed};
and by taking $\W$ to be an ordinary monoidal category, viewed as a locally
discrete monoidal bicategory, whereupon $\W$-bicategories reduce to categories
enriched in $\W$. An account of the general theory of enriched bicategories can
be found in \cite{Lack}, but we will need sufficiently little of it that we can
easily arrange for our account to be self-contained:
\begin{Defn}
A \defn{locally cubical bicategory} is a bicategory enriched in the monoidal
2-category $\cat{DblCat}$. Explicitly, it is given by the following data:
\begin{enumerate}[(LDD1)]
\item[(LDD1)] A collection $\ob {\mathfrak B}$ of objects;
\item[(LDD2)] For every pair $A, B \in \ob {\mathfrak B}$, a pseudo double
    category ${\mathfrak B}(A, B)$;
\item[(LDD3)] For every $A \in \ob {\mathfrak B}$, a unit homomorphism
\[\elt{I_x} \colon 1 \to {\mathfrak B}(A, A)\text;\]
\item[(LDD4)] For every triple $A, B, C \in \ob {\mathfrak B}$, a
    composition homomorphism
\[\otimes \colon {\mathfrak B}(B, C) \times {\mathfrak B}(A, B) \to {\mathfrak B}(A, C)\text;\]
\item[(LDD5)] For every pair $A, B \in \ob {\mathfrak B}$, invertible
    vertical transformations
\[\cd[@C+1em]{
{\mathfrak B}(A, B) \times {\mathfrak B}(A, A) \ar[dr]_{\otimes} & {\mathfrak B}(A, B) \rtwocell[0.36]{dl}{r}
\ltwocello[0.36]{dr}{l} \ar[r]^-{\elt{I_B} \times 1} \ar[l]_-{1 \times \elt{I_A}} \ar[d]_{1} &
{\mathfrak B}(B, B) \times {\mathfrak B}(A, B)\text; \ar[dl]^{\otimes} \\ & {\mathfrak B}(A, B) & {}
 }\]
\item[(LDD6)] For every quadruple $A, B, C, D \in \ob {\mathfrak B}$, an
    invertible vertical transformation
\[\cd[@C+1em]{
{\mathfrak B}(C, D) \times {\mathfrak B}(B, C) \times {\mathfrak B}(A, B) \ar[d]_{\otimes \times 1} \ar[r]^-{1 \times \otimes}
\dtwocell{dr}{a} &
{\mathfrak B}(C, D) \times {\mathfrak B}(A, C) \ar[d]^{\otimes}\\
{\mathfrak B}(B, D) \times {\mathfrak B}(A, B) \ar[r]_-{\otimes} & {\mathfrak B}(A, D)\text.}\]
\end{enumerate}

Subject to the following two axioms:

\begin{enumerate}[(LDA1)]
\item[(LDA1)] For each triple of objects $A, B, C$ of ${\mathfrak B}$, the
    following pasting equality holds:
\[\cd[@+1em@C+1em]{
  {\mathfrak B}^2 \ar[r]^{1 \times \elt{I_B} \times 1} \ar[d]_{1} \dtwocell{dr}{1 \times l} & {\mathfrak B}^3 \ar[r]^{1 \times \otimes} \ar[d]^{1 \times \otimes} & {\mathfrak B}^2 \ar[d]^{\otimes}\\
  {\mathfrak B}^2 \ar[r]_1 & {\mathfrak B}^2 \ar[r]_{\otimes} & {\mathfrak B}
} \quad = \quad \cd[@+1em@C+1em]{
  {\mathfrak B}^2 \ar[r]^{1 \times \elt{I_B} \times 1} \ar[d]_{1} \dtwocell{dr}{r \times 1} & {\mathfrak B}^3 \ar[r]^{1 \times \otimes} \ar[d]^{\otimes \times 1} \dtwocell{dr}{a} & {\mathfrak B}^2 \ar[d]^{\otimes}\\
  {\mathfrak B}^2 \ar[r]_1 & {\mathfrak B}^2 \ar[r]_{\otimes} & {\mathfrak
B}\text, }\] where $\mathfrak B^2$ and $\mathfrak B^3$ abbreviate the appropriate products of hom-double categories;
\item[(LDA2)] For each quintuple of objects $A, B, C, D, E$ of ${\mathfrak
    B}$, the following pasting equality holds:
\[
\cd[@!@+1em]{
  {\mathfrak B}^4 \ar[r]^{1 \times 1 \times \otimes} \ar[d]_{\otimes\times  1\times  1} \twocong{dr} & {\mathfrak B}^3 \ar[dr]^{1 \times \otimes} \ar[d]^{\otimes \times 1}\\
  {\mathfrak B}^3 \ar[r]_{1 \times \otimes} \ar[dr]_{\otimes \times 1}  & {\mathfrak B}^2 \ar[dr]^{\otimes} \drtwocell{d}{a} \dtwocell{r}{a} & {\mathfrak B}^2 \ar[d]^{\otimes} \\
  & {\mathfrak B}^2 \ar[r]_{\otimes} & {\mathfrak B}
} \quad = \quad \cd[@!@+1em]{
  {\mathfrak B}^4 \ar[r]^{1 \times 1 \times \otimes} \ar[dr]|{1 \times \otimes \times 1} \ar[d]_{\otimes\times  1\times  1} & {\mathfrak B}^3 \ar[dr]^{1 \times \otimes} \drtwocell{d}{1 \times a}\\
  {\mathfrak B}^3 \ar[dr]_{\otimes \times 1} \drtwocell[0.4]{r}{a \times 1} & {\mathfrak B}^3 \dtwocell{dr}{a} \ar[d]^{\otimes \times 1} \ar[r]_{1 \times \otimes} & {\mathfrak B}^2 \ar[d]^{\otimes} \\
  & {\mathfrak B}^2 \ar[r]_{\otimes} & {\mathfrak B}\text,
}
\]
where we observe the same convention regarding $\mathfrak B^4$, $\mathfrak
B^3$ and $\mathfrak B^2$.
\end{enumerate}
\end{Defn}

It may be helpful to extract a description of the various sorts of composition
that a $\cat{DblCat}$-bicategory possesses. The 0-cells, 1-cells and vertical
2-cells form an ordinary bicategory. Next come the the horizontal 2-cells,
which can be composed with each other along either a 1-cell boundary or a
0-cell boundary, with both compositions being associative up to an invertible
globular 3-cell; moreover, the corresponding ``middle four interchange'' law
only holds up to an invertible globular 3-cell. Finally, the 3-cells themselves
can be composed with each other along the two different types of 2-cell
boundary and along 0-cell boundaries; and these operations are strictly
associative modulo the associativity of the boundaries.

A \emph{one-object} locally cubical bicategory amounts to a
\defn{monoidal double category} \cite{doubleclubs,GP2,mike}~--~that is, a pseudo double category with an
up-to-isomorphism tensor product on it. In particular, any double category with \emph{finite
products} in the appropriate double categorical sense\footnote{By which we mean a
\emph{pseudo-functorial choice of double products} in the sense of \cite{GP1}. Such pseudo
double categories are slightly stricter versions of the \emph{cartesian bicategories} of
\cite{cb2}, which, although they are presented in a globular way, are essentially cubical
structures.}\ \ becomes a monoidal double category under the cartesian tensor product. The
double categories $\mathfrak{Cat}$, $\mathfrak{Span}(\C)$ (where $\C$ is a category with
finite limits) and $\mathfrak{Rng}$ are all monoidal in this way: though in the case of
$\mathfrak{Rng}$, there is another natural monoidal structure which is derived from the tensor
product on the category of rings.

For a non-degenerate example of a locally cubical bicategory, we turn to
$\cat{DblCat}$ itself. As demonstrated in \cite{doubleclubs}, we may define an
internal hom 2-functor
\[\HOM{\, \thg , ?} \colon \cat{DblCat}^\op \times
\cat{DblCat} \to \cat{DblCat}\label{hom}\] for which $\HOM{\mathfrak C, \mathfrak D}$ is the
following double category. Its objects are homomorphisms $\mathfrak C \to \mathfrak D$, and
its vertical 1-cells $\alpha \colon F \Rightarrow G$ are the vertical transformations between
them. Its horizontal 1-cells $\alpha \colon F \Tor G$ are the \emph{horizontal pseudo-natural
transformations}, whose components at an object $x \in \mathfrak C$ are given by horizontal
1-cells $\alpha_x \colon Fx \tor Gx$ of $\mathfrak D$, together with pseudo-naturality data
like that for a pseudo-natural transformation of bicategories; and indeed, in the case that
$\mathfrak C$ and $\mathfrak D$ are bicategories the two notions coincide. Finally, the
2-cells of $\HOM{\mathfrak C, \mathfrak D}$ are the \emph{cubical modifications}, which are
bounded by two horizontal and two vertical transformations and whose basic data consists of
giving, for each object of the source, a 2-cell of the target bounded by the components of
these transformations: Definition \ref{cubbicat} below makes this explicit in the special case
where $\mathfrak C$ and $\mathfrak D$ are bicategories.

When we say that $\HOM{\thg, ?}$ acts as an internal hom, we are affirming a
universal property: namely, that for each $\mathfrak C$ the 2-functor $(\thg)
\times \mathfrak C \colon \cat{DblCat} \to \cat{DblCat}$ is left biadjoint to
$\HOM{\, \mathfrak C, \thg}$, so that what we have is a \emph{biclosed}
monoidal bicategory in the sense of \cite{monoidalbicats}. Now, in \cite{Lack},
it is demonstrated that, just as any closed monoidal category can be viewed as
a category enriched over itself, so any biclosed monoidal bicategory can be
viewed as a bicategory enriched over itself, with the hom-objects being given
by the biclosed structure. Applying this result to the monoidal 2-category
$\cat{DblCat}$, we obtain a locally cubical bicategory $\mathfrak{DblCat}$,
with 0-cells being the pseudo double categories; 1-cells, the homomorphisms;
vertical 2-cells, the vertical transformations; horizontal 2-cells, the
horizontal pseudo-natural transformations; and 3-cells the cubical
modifications.

In particular, if we restrict our attention to those pseudo double categories
lying in the image of the embedding $\mathbb U \colon \cat{Bicat} \to
\cat{DblCat}$ then we obtain:
\begin{Cor}\label{bicategories-as-dblcatbicategories}
There is a locally cubical bicategory $\mathfrak{Bicat}$ which has as 0-cells,
bicategories; as 1-cells, homomorphisms; as vertical 2-cells, bicategorical
icons; as horizontal 2-cells, pseudo-natural transformations; and as 3-cells,
cubical modifications.
\end{Cor}
Whilst the 0-, 1-, and 2-cells of $\mathfrak{Bicat}$ are familiar, the same is not true of the
3-cells; and since we will need them in Definition~\ref{cubicalpseudo} below, it is worth
giving an explicit description.
\begin{Defn}\label{cubbicat}
Suppose that $F, G, H, K \colon \B \to \C$ are homomorphisms of bicategories;
that \mbox{$\alpha \colon F \Tor G$} and $\beta \colon H \Tor K$ are
pseudo-natural transformations; and that \mbox{$\gamma \colon F \Rightarrow H$}
and $\delta \colon G \Rightarrow K$ are bicategorical icons. Then a
\defn{cubical modification}
\[\cd{
F \ar@2[d]_{\gamma} \ar@2[r]|-{\object@{|}}^{\alpha} \dthreecell{dr}{\Gamma}& G \ar@2[d]^{\delta}  \\
H \ar@2[r]|-{\object@{|}}_{\beta} & K}\] is given by specifying, for every
object $A \in \B$ a 2-cell $\Gamma_A \colon \alpha_A \Rightarrow \beta_A$, in
such a way that for every 1-cell $f \colon A \to B$ of $\B$, the following
pasting equality holds:
\[ \cd[@!]{
    & FB
      \ar[r]^{\alpha_B}  \rtwocell{d}{\alpha_f} \ar@{<-}[dl]_{Ff}
    & GB
      \ar@{=}[d] \ar@{<-}[dl]|{Gf} \\
    FA
      \ar[r]^{\alpha_A} \ar@{=}[d] \dtwocell{dr}{\Gamma_A}
    & GA
      \ar@{=}[d] \dtwocell{r}{\delta_f}
    & KB
      \ar@{<-}[dl]^{Kf}
    \\ HA
      \ar[r]_{\beta_A}
    &  KA\text;
} \quad = \quad \cd[@!]{
    & FB
      \ar[r]^{\alpha_B}  \ar@{=}[d] \dtwocell{dr}{\Gamma_B} \ar@{<-}[dl]_{Ff}
    & GB
      \ar@{=}[d] \\
    FA
       \ar@{=}[d] \dtwocell{r}{\gamma_f}
    & HB
      \ar[r]_{\beta_B} \ar@{<-}[dl]|{Hf} \rtwocell{d}{\beta_f}
    & KB
      \ar@{<-}[dl]^{Kf}
    \\ HA
      \ar[r]_{\beta_A}
    &  KA\text.
}\]
\end{Defn}
In particular, to give a \emph{globular} 3-cell of $\mathfrak{Bicat}$ is
precisely to give a modification between pseudo-natural transformations in the
standard sense; and so $\mathfrak{Bicat}$ is rich enough to encode faithfully
all the cells and all of the forms of composition which feature in the
tricategory of bicategories, but is able to do so using coherence whose
complexity does not rise above the bicategorical level.

Pleasing as this is, we should note that not every tricategory can be reduced
to a locally cubical bicategory in this way; for example, given a bicategory
$\B$ with bipullbacks, we may form the tricategory $\cat{Span}(\B)$ of spans in
$\B$. In this tricategory, 1-cell composition is given by bipullback, and so is
only determined up-to-equivalence, rather than up-to-isomorphism; so evidently,
it will be inexpressible as a locally cubical bicategory.

\begin{Rk}
There are two canonical ways of forming a tricategory of bicategories,
corresponding to the two canonical ways of composing a pair of strong
transformations along a 0-cell boundary: however, Proposition
\ref{bicategories-as-dblcatbicategories} exhibited a single canonical locally
cubical bicategory of bicategories. The discrepancy is resolved if we observe
that to obtain this $\cat{DblCat}$-bicategory we must fix a choice of biclosed
structure on $\cat{DblCat}$, and that there are two canonical ways of doing
this, depending on how we choose the counit maps $\HOM{\mathfrak B, \mathfrak
C} \times \mathfrak B \to \mathfrak C$ for the biadjunctions in question.
\end{Rk}

\section{A locally cubical bicategory of tricategories}\label{Sec:local} We now
return to our study of tricategories with the goal of forming them into a
locally cubical bicategory. The result we will prove in this section is:

\begin{Thm}\label{loccubtricat}
There is a locally cubical bicategory $\mathfrak{Tricat}_3$ with 0-cells being
tricategories; 1-cells, lax homomorphisms; vertical 2-cells, ico-icons;
horizontal 2-cells, pseudo-icons; and 3-cells, cubical pseudo-icon
modifications.
\end{Thm}

We have already met the lax homomorphisms (Definition \ref{laxhomom}), the
ico-icons (Definition \ref{icoicon}) and the pseudo-icons (Definition
\ref{pseudoicon}); however, we have not yet introduced the \emph{cubical
pseudo-icon modifications}. These generalise the (globular) pseudo-icon
modifications of Definition \ref{pseudoiconmod} as follows:

\begin{Defn}\label{cubicalpseudo}
Let $F$, $G$, $F'$, $G' \colon \S \to \T$ be lax homomorphisms of
tricategories, let $\alpha \colon F \Tor G$ and $\beta \colon F' \Tor G'$ be
pseudo-icons, and let $\gamma \colon F \Rightarrow F'$ and $\delta \colon G
\Rightarrow G'$ be ico-icons. Then a
\defn{cubical pseudo-icon modification}\vskip-\baselineskip
\[\cd{
F \ar@2[d]_{\gamma} \ar@2[r]|-{\object@{|}}^{\alpha} \dthreecell{dr}{\Gamma}& G \ar@2[d]^{\delta}  \\
F' \ar@2[r]|-{\object@{|}}_{\beta} & G'}\] consists in the following data:
\begin{enumerate}[(MD1)]
\item[(MD1)] For each $A, B$ in $\S$, a cubical modification (cf.
    Definition \ref{cubbicat})\vskip-\baselineskip
\[\cd{
F_{A, B} \ar@2[d]_{\gamma_{A, B}} \ar@2[r]|-{\object@{|}}^{\alpha_{A, B}} \dthreecell{dr}{\Gamma_{A, B}}& G_{A, B} \ar@2[d]^{\delta_{A, B}}  \\
F'_{A, B} \ar@2[r]|-{\object@{|}}_{\beta_{A, B}} & G'_{A, B}}\] (and so in
 particular, for each 1-cell $f \colon A \to B$ of $\S$, a
3-cell $\Gamma_f \colon \alpha_f \Rrightarrow \beta_f$ of $\T$);
\end{enumerate}
subject to the following axioms:
\begin{enumerate}[(MA1)]
\item[(MA1)] For each object $A$ of $\S$, the following pasting equality
    holds:
\[\cd[@u@!@+1em]{
    I_{F'A} \ar@{=}[r] \ar@{=}[d] \ar@2[dr]_{\iota^{F'}_A} &
    I_{FA} \ar@2[dr]^{\iota^F_A} \dthreecell[0.45]{d}{M^{\gamma}_A} \\
    I_{G'A} \ar@2[dr]_{\iota^{G'}_A} \dlthreecell{r}{M^{\beta}_A}  &
    F'I_A \ar@2[d]^{\beta_{I_A}} \ar@{=}[r] \dthreecell{dr}{\Gamma_{I_A}} &
    FI_A \ar@2[d]^{\alpha_{I_A}}\\
    & G'I_A \ar@{=}[r] & GI_A
}\quad = \quad \cd[@u@!@+1em]{
    I_{F'A} \ar@{=}[r] \ar@{=}[d] &
    I_{FA} \ar@{=}[d] \ar@2[dr]^{\iota^F_A} \\
    I_{G'A} \ar@{=}[r] \ar@2[dr]_{\iota^{G'}_A}  &
    I_{GA} \ar@2[dr]_{\iota^G_A} \dlthreecell{r}{M^{\alpha}_A} \dthreecell[0.45]{d}{M^{\delta}_A} &
    FI_A \ar@2[d]^{\alpha_{I_A}}\\
    & G'I_A \ar@{=}[r] & GI_A\text;
}
\]
\item[(MA2)] For each pair of composable 1-cells $f \colon A \to B$, $g
    \colon B \to C$ of $\S$, the following pasting equality holds:
\end{enumerate}
\[\cd[@u@!@+1.5em@R-1em]{
    F'g . F'f \ar@{=}[r] \ar@2[d]_{\beta_g . \beta_f} \ar@2[dr]_{\chi^{F'}_{f, g}} &
    Fg . Ff \ar@2[dr]^{\chi^F_{f, g}} \dthreecell{d}{\Pi^{\gamma}_{f, g}} \\
    G'g . G'f \ar@2[dr]_{\chi^{G'}_{f, g}} \dlthreecell{r}{\Pi^{\beta}_{f, g}}  &
    F'(gf) \ar@2[d]^{\beta_{gf}} \ar@{=}[r] \dthreecell{dr}{\Gamma_{gf}} &
    F(gf) \ar@2[d]^{\alpha_{gf}}\\
    & G'(gf) \ar@{=}[r] & G(gf)
}\quad = \quad \cd[@u@!@+1.5em@R-1em]{
    F'g . F'f \ar@{=}[r] \ar@2[d]_{\beta_g . \beta_f} \dthreecell{dr}{\Gamma_g . \Gamma_f} &
    Fg . Ff \ar@2[d]^{\alpha_g . \alpha_f} \ar@2[dr]^{\chi^F_{f, g}} \\
    G'g . G'f \ar@{=}[r] \ar@2[dr]_{\chi^{G'}_{f, g}}  &
    Gg . Gf \ar@2[dr]_{\chi^G_{f, g}} \dlthreecell{r}{\Pi^{\alpha}_{f, g}} \dthreecell[0.45]{d}{\Pi^{\delta}_{f, g}} &
    F(gf) \ar@2[d]^{\alpha_{gf}}\\
    & G'(gf) \ar@{=}[r] & G(gf)\text.
}
\]
\end{Defn}
\noindent The first step in the proof of Theorem \ref{loccubtricat} will be to
give the local structure:
\begin{Prop}\label{localstruct}
Let $\S$ and $\T$ be tricategories. Then the lax homomorphisms, ico-icons,
pseudo-icons and cubical pseudo-icon modifications from $\S$ to $\T$ form a
pseudo double category $\mathfrak{Tricat}_3(\S, \T)$.
\end{Prop}
\begin{proof}
Underlying each lax homomorphism, ico-icon, pseudo-icon or cubical pseudo-icon
modification is an indexed family of homomorphisms of bicategories,
bicategorical icons, pseudo-natural transformations, or cubical modifications,
respectively: thus our approach will be to lift the compositional structure
from the pseudo double categories $\HOM{\mathfrak C, \mathfrak D}$ as defined
preceding Corollary~\ref{bicategories-as-dblcatbicategories}.

We begin with the vertical structure of $\mathfrak{Tricat}_3(\S, \T)$. We have already seen in
Section~\ref{bicat-tricat} that the lax homomorphisms and ico-icons from $\S$ to $\T$ form a
category; we must show the same is true of the pseudo-icons and the cubical pseudo-icon
modifications. So for each pseudo-icon $\alpha \colon F \Tor G$ we must give a cubical
pseudo-icon modification
\[\cd{
  F  \ar@2[d]_{\id_F} \ar@2[r]|-{\object@{|}}^{\alpha} \dthreecell{dr}{\id_{\alpha}}& G  \ar@2[d]^{\id_G}  \\
  F \ar@2[r]|-{\object@{|}}_{\alpha} & G\text;}\] which we take to be given by the identity family of cubical
modifications $(\id_\alpha)_{A, B} = \id_{\alpha_{A, B}} \colon \alpha_{A, B} \Rrightarrow
\alpha_{A, B}$. The axioms (MA1) and (MA2) are clear, since every occurence of $\Gamma$
reduces to an identity 3-cell. Next, given cubical pseudo-icon modifications
\[\cd{
  F  \ar@2[d]_{\sigma} \ar@2[r]|-{\object@{|}}^{\alpha} \dthreecell{dr}{\Gamma}& F'  \ar@2[d]^{\sigma'}  \\
  G  \ar@2[d]_{\tau} \ar@2[r]|-{\object@{|}}^{\beta} \dthreecell{dr}{\Delta}& G'  \ar@2[d]^{\tau'}  \\
  H \ar@2[r]|-{\object@{|}}_{\gamma} & H'\text;}\]
we must provide a vertical composite $\Delta\Gamma \colon \alpha \Rrightarrow \gamma$, which
we do by composing their underlying families of cubical modifications:
\[(\Delta\Gamma)_{A, B} = \Delta_{A, B}. \Gamma_{A, B} \colon \alpha_{A, B} \Rrightarrow \gamma_{A, B}\text.\]
Now the axioms (MA1) and (MA2) follow from an application of the corresponding
axiom for $\Delta$ followed by the corresponding axiom for $\Gamma$.
Associativity and unitality of this composition follow from that for
composition of cubical modifications.

We next describe the horizontal identities of $\mathfrak{Tricat}_3(\S, \T)$. Firstly, for each
lax homomorphism $F \colon \S \to \T$, we must give an identity pseudo-icon $\DI F \colon F
\Tor F$. This has (ID1) given by the family $(\DI F)_{A, B} = \id_{F_{A, B}} \colon F_{A, B}
\Rightarrow F_{A, B}$ whilst $M^{\DI F}_A$ and $\Pi^{\DI F}_{A, B, C}$ are given by unnamed
coherence isomorphisms in the hom-bicategories of $\T$. Secondly, for each ico-icon $\alpha
\colon F \Rightarrow G$, we must give a cubical pseudo-icon modification
\[\cd{
  F  \ar@2[d]_{\alpha} \ar@2[r]|-{\object@{|}}^{\DI F} \dthreecell{dr}{\DI \alpha}& F  \ar@2[d]^{\alpha}  \\
  G \ar@2[r]|-{\object@{|}}_{\DI G} & G\text;}\] which we do by taking (MD1) to be given by the identity family of cubical
modifications $\id_{\id_{Ff}} \colon \id_{Ff} \Rrightarrow \id_{Gf}$. Each of the axioms
(IA1)--(IA3) for $\DI F$ and (MA1)--(MA2) for $\DI \alpha$ now asserts that some 3-cell is
equal to itself when pasted with such unnamed coherence cells, and this follows from coherence
for bicategories. Finally, we must check functoriality of $\DI{(\thg)}$, which is immediate.

We now come to the horizontal composition of $\mathfrak{Tricat}_3(\S, \T)$. First, for each
pair of pseudo-icons $\alpha \colon F \Tor G$ and $\beta \colon G \Tor H$, we must give a
composite pseudo-icon $\beta \alpha \colon F \Tor H$. We do this as follows:
\begin{enumerate}[(ID1)]
\item[(ID1)] $(\beta\alpha)_{A, B} = \beta_{A, B}.\alpha_{A, B} \colon
    F_{A, B} \Rightarrow H_{A, B}$;
\item[(ID2)] $M^{\beta\alpha}_A$ is the pasting:
\[
\cd[@C+1em]{
    I_{FA} \ar@2[r]^{\iota^F_A} \ar@{=}[d] \dthreecell{dr}{M^{\alpha}_A} & FI_A \ar@2[d]^{\alpha_{I_A}} \\
    I_{GA} \ar@2[r]|{\iota^G_A} \ar@{=}[d] \dthreecell{dr}{M^{\beta}_A} & GI_A \ar@2[d]^{\beta_{I_A}} \\
    I_{HA} \ar@2[r]_{\iota^H_A} & HI_A\text;
}\]
\item[(ID3)] $\Pi^\alpha_{A, B, C}$ is the pasting:
\[
\cd{
  F(\thg) \otimes F(?) \ar@2[r]^{\chi^F} \ar@2[d]|{\alpha_{(\thg)}\otimes \alpha_{(?)}} \dthreecell{dr}{\Pi_{A, B, C}^\alpha} \ar@/_5em/@2[dd]_{(\beta\alpha)_{(\thg)} \otimes (\beta\alpha)_{(?)}} &
  F\big((\thg) \otimes (?)\big) \ar@2[d]^{\alpha_{((\thg) \otimes (?))}} \\
  G(\thg)\otimes G(?) \twocong[-0.4]{r} \ar@2[r]^{\chi^G} \ar@2[d]|{\beta_{(\thg)}\otimes \beta_{(?)}} \dthreecell{dr}{\Pi_{A, B, C}^\beta} &
  G\big((\thg) \otimes (?)\big) \ar@2[d]^{\beta_{((\thg) \otimes (?))}} \\
  H(\thg)\otimes H(?) \ar@2[r]_{\chi^H} &
  H\big((\thg) \otimes (?)\big)\text.
}
\]
\end{enumerate}
Showing that these data satisfy axioms (IA1)--(IA3) is almost as simple as
placing the corresponding diagrams for $\beta$ and $\alpha$ alongside each
other; though not quite, since there are a number of auxiliary coherence
results we need to prove first. For instance, in order to prove (IA1) we must
show that:
\begin{multline*}
\cd{
  I_{FB} . Ff \ar@2[r]^{\iota^F .1} \ar@2[d]|{1 . \alpha_f} \dthreecell{dr}{\overline{M^\alpha . 1}} \ar@/_5em/@2[dd]_{1 . (\beta\alpha)_f} &
  FI_B . Ff \ar@2[d]^{\alpha_{I_B} . \alpha_f} \ar@/^5em/@2[dd]^{(\beta\alpha)_{I_B} . (\beta\alpha)_f}\\
  I_{GB} . Gf \twocong[-0.4]{r} \ar@2[r]^{\iota^G .1} \ar@2[d]|{1 . \beta_f} \dthreecell{dr}{\overline{M^\beta . 1}} &
  GI_B . Gf \twocong[-0.4]{l} \ar@2[d]^{\beta_{I_B} . \beta_f} \\
  I_{HB} . Hf \ar@2[r]_{\iota^H . 1} &
  HI_B . Hf
} \\ = \quad \cd[@+1em]{
  I_{FB} . Ff \ar@2[r]^{\iota^F .1} \ar@2[d]_{1 . (\beta\alpha)_f} \dthreecell{dr}{\overline{M^{\beta\alpha} .
  1}} &
  FI_B . Ff \ar@2[d]^{(\beta\alpha)_{I_B} . (\beta\alpha)_f} \\
  I_{HB} . Hf \ar@2[r]_{\iota^H .1} &
  HI_B . Hf\text;
}
\end{multline*}
holds; and similarly for (IA2) and (IA3). These derivations are straightforward
bicategorical manipulations and left to the reader.

Secondly, for each diagram of cubical pseudo-icon modifications
\[\cd{
F \ar@2[d]_{\sigma} \ar@2[r]|-{\object@{|}}^{\alpha} \dthreecell{dr}{\Gamma}& G\ar@2[r]|-{\object@{|}}^{\beta} \dthreecell{dr}{\Delta} \ar@2[d]^{\tau} & H \ar@2[d]^{\upsilon}  \\
F' \ar@2[r]|-{\object@{|}}_{\alpha'} & G' \ar@2[r]|-{\object@{|}}_{\beta'} & H'}\] we must
give a cubical pseudo-icon modification $\Delta \ast \Gamma \colon \beta \alpha \Rrightarrow
\beta' \alpha' \colon F \Rightarrow H$. We do this by taking
\[(\Delta \ast \Gamma)_{A, B} = \Delta_{A, B} \ast \Gamma_{A, B} \colon
\beta_{A, B} \alpha_{A, B} \Rrightarrow \beta'_{A, B}. \alpha'_{A, B}\text,\] where $\ast$ on
the right-hand side is horizontal composition of cubical modifications in the pseudo double
category $\HOM{\S(A, B), \T(FA, FB)}$. Explicitly, for any 1-cell $f$ of the tricategory $\S$,
the 3-cell $(\Delta \ast \Gamma)_f$ of the tricategory $\T$ is given by the pasting
\[\cd[@u]{
    Ff \ar@2[d]_{\alpha'_f} \ar@{=}[r] \dthreecell{dr}{\Gamma_f} & Ff \ar@2[d]^{\alpha_f}\\
    Gf \ar@2[d]_{\beta'_f} \ar@{=}[r] \dthreecell{dr}{\Delta_f} & Gf \ar@2[d]^{\beta_f}\\
    Hf \ar@{=}[r] & Hf\text,
}\] and thus (MA1) and (MA2) for $\Delta \ast \Gamma$ follow by placing the corresponding
axioms for $\Gamma$ and $\Delta$ beside each other, together with some very simple
manipulation with unnamed coherence cells. Finally, we must check functoriality of the
horizontal composition functor, which is just the middle-four interchange law. This will hold
in $\mathfrak{Tricat}_3(\S, \T)$ because it does in each double category $\HOM{\S(A, B),
\T(FA, FB)}$.

It remains only to give the unitality and associativity constraints for the pseudo double
category $\mathfrak{Tricat}_3(\S, \T)$. So let there be given pseudo-icons $\alpha \colon F
\Rightarrow G$, $\beta \colon G \Rightarrow H$ and $\gamma \colon H \Rightarrow K$. Then:
\begin{itemize}
\item The associativity constraint $\a_{\alpha, \beta, \gamma} \colon
    (\gamma\beta)\alpha \Rrightarrow \gamma(\beta\alpha)$ has  component
    modification $(\a_{\alpha, \beta, \gamma})_{A, B}$ given by the
    associativity constraint $\a_{\alpha_{A, B}, \beta_{A, B}, \gamma_{A,
    B}}$ in the double category $\HOM{\S(A, B), \T(FA, FB)}$;
\item The left unitality constraint $\l_\alpha \colon \id_G . \alpha
    \Rrightarrow \alpha$ has component modification $(\l_\alpha)_{A, B}$
    given by the left unitality constraint $\l_{\alpha_{A, B}}$ in
    $\HOM{\S(A, B), \T(FA, FB)}$;
\item The right unitality constraint $\r_\alpha \colon \alpha . \id_G
    \Rrightarrow \alpha$ has component modification $(\r_\alpha)_{A, B}$
    given by the right unitality constraint $\r_{\alpha_{A, B}}$ in
    $\HOM{\S(A, B), \T(FA, FB)}$.
\end{itemize}
The naturality of these constraints in $\alpha$, $\beta$ and $\gamma$ is
inherited from the hom-double categories $\HOM{\S(A, B), \T(FA, FB)}$; and that
these data satisfy the axioms (MA1) and (MA2) is also straightforward. In the
case of $\a_{\alpha, \beta, \gamma}$, for example, we see that $M^{\gamma(
\beta\alpha)}_A$ and $\Pi^{\gamma(\beta\alpha)}_{f, g}$ can be obtained from
$M^{(\gamma \beta) \alpha}_A$ and $\Pi^{(\gamma\beta)\alpha}_{f, g}$ by pasting
with unnamed coherence isomorphisms; but the components of $\a_{\alpha, \beta,
\gamma}$ are built from the selfsame coherence isomorphisms, and so the result
follows from the coherence theorem for bicategories.
\end{proof}
In order for the pseudo double categories $\mathfrak{Tricat}_3(\S, \T)$ we have
just defined to provide homs for the locally cubical bicategory
$\mathfrak{Tricat}_3$, we must define double homomorphisms which provide
top-level identities and composition. The double homomorphism
\[\elt{I_\T} \colon 1 \to \mathfrak{Tricat}_3(\T, \T)\] is straightforward; it sends the unique object of the terminal pseudo double category
to the identity lax homomorphism $\id_\T \colon \T \to \T$, the unique vertical 1-cell to the
identity ico-icon on $\id_T$; the unique horizontal 1-cell to the identity pseudo-icon on
$\id_\T$; and the unique 2-cell to the identity cubical pseudo-icon modification on this. The
coherence data for this homomorphism is obtained from unitality constraints in
$\mathfrak{Tricat}_3(\T, \T)$, and so the homomorphism axioms follow from coherence for
bicategories. We must now give the composition double homomorphism
\[\otimes \colon \mathfrak{Tricat}_3(\T, \U) \times \mathfrak{Tricat}_3(\S, \T) \to \mathfrak{Tricat}_3(\S, \U)\text.\]
The general approach to defining this will be similar to that adopted in the proof of
Proposition \ref{localstruct}. There, we defined the compositional structure on
$\mathfrak{Tricat}_3(\S, \T)$ by lifting it from the pseudo double categories $\HOM{\S(A, B),
\T(FA, FB)}$: here, we will define $\otimes$ by lifting the double homomorphisms
\begin{multline*}
    \HOM{\T(FA, FB), \U(GFA, GFB)} \times \HOM{\S(A, B), \T(FA,
FB)} \\ \longrightarrow \ \HOM{\S(A, B), \U(GFA, GFB)}
\end{multline*}
which provide composition in the locally cubical bicategory $\mathfrak{Bicat}$ of Proposition
\ref{bicategories-as-dblcatbicategories}.

In detail, $\otimes$ is given as follows. On objects and vertical 1-cells, it
is given by the composition law for $\cat{Tricat}_2$. On horizontal 1-cells, we
consider pseudo-icons \mbox{$\alpha \colon F \Tor F' \colon \S \to \T$} and
$\beta \colon G \Tor G' \colon \T \to \U$, for which the composite pseudo-icon
$\beta \otimes \alpha \colon GF \Tor G'F'$ is given as follows:
\begin{enumerate}[(ID1)]
\item[(ID1)] $(\beta \otimes \alpha)_{A, B} = \beta_{FA, FB} \otimes
    \alpha_{A, B}$, where $\otimes$ on the right-hand side is one of the
    two canonical choices for horizontal composition of pseudo-natural
    transformations; for concreteness let us take
\[(\beta \otimes \alpha)_f = \cd{GFf \ar@2[r]^{\beta_{Ff}} & G'Ff \ar@2[r]^{G'\alpha_f} & G'F'f}\]
and
\[(\beta \otimes \alpha)_\theta \quad = \quad \cd{
    GFf \ar@2[r]^{GF \theta} \ar@2[d]_{\beta_{Ff}} \dthreecell{dr}{\beta_{F\theta}}& GFg \ar@2[d]^{\beta_{Fg}} \\
    G'Ff \ar@2[r]^{G'F \theta} \ar@2[d]_{G'\alpha_f} \dthreecell{dr}{G'\alpha_{\theta}}& G'Fg \ar@2[d]^{G'\alpha_g} \\
    G'F'f \ar@2[r]_{G'F' \theta} & G'F'g\text{;}
}\]
\item[(ID2)] $M^{\beta \otimes \alpha}_A$ is the following 3-cell:
\[\cd{
    I_{GFA}
      \ar@2[r]^{\iota^G} \ar@{=}[d] \dthreecell{dr}{M^{\beta}_{FA}}
    & GI_{FA}
      \ar@2[d]^{\beta_{I_{FA}}} \ar@2[r]^{G\iota^F} \dthreecell{dr}{\beta_{\iota^F}}
    & GFI_A
      \ar@2[d]^{\beta_{FI_A}}
    \\ I_{G'FA}
      \ar@2[r]^{\iota^{G'}} \ar@{=}[d] \twoeq{dr}
    & G'I_{FA}
      \ar@{=}[d] \ar@2[r]^{G'\iota^F} \dthreecell{dr}{\overline{G'M^{\alpha}_A}}
    & G'FI_A
      \ar@2[d]^{G'\alpha_{I_A}}
    \\ I_{G'F'A}
      \ar@2[r]_{\iota^{G'}}
    & G'I_{F'A}
      \ar@2[r]_{G'\iota^{F'}}
    & G'F'I_A\text;}
\]
\item[(ID3)] $\Pi^{\beta \otimes \alpha}_{A, B, C}$ is the pseudo-natural
    transformation with the following components:
\[
\cd[@C+1em]{
    GFg . GFf
      \ar@/_5em/@2[dd]_{(\beta \otimes \alpha)_g . (\beta \otimes \alpha)_f} \ar@2[r]^{\chi^G} \ar@2[d]|{\beta_{Fg} . \beta_{Ff}} \dthreecell{dr}{\Pi^{\beta}_{Ff, Fg}}
    & G(Fg . Ff)
      \ar@2[r]^{G\chi^F} \ar@2[d]^{\beta_{Fg . Ff}} \dthreecell{dr}{\beta_{\chi^F}}
    & GF(gf)
      \ar@2[d]^{\beta_{F(gf)}} \\
    G'Fg . G'Ff \twocong[-0.4]{r}
      \ar@2[r]^{\chi^{G'}} \ar@2[d]|{G'\alpha_g . G'\alpha_f} \twocong{dr}
    & G'(Fg . Ff)
      \ar@2[r]^{G'\chi^F} \ar@2[d]|{G'(\alpha_g . \alpha_f)} \dthreecell{dr}{\overline{G'\Pi^\alpha_{f, g}}}
    & G'F(gf)
      \ar@2[d]^{G'\alpha_{gf}} \\
    G'F'g . G'F'f
      \ar@2[r]_{\chi^{G'}}
    & G'(F'g . F'f)
      \ar@2[r]_{G'\chi^{F'}}
    & G'F'(gf)\text.
}\]
\end{enumerate}
The proof that these data satisfy axioms (IA1)--(IA3) consists once again in
building large cubes or hexagonal prisms from smaller ones, together with some
simple manipulation with unnamed coherence cells: and once again, we leave this
task to the reader.

Finally we must define the action of $\otimes$ on pseudo-icon modifications.
Given two such:
\[
\cd{
 F \ar@2[d]_{\sigma} \ar@2[r]|-{\object@{|}}^{\alpha} \dthreecell{dr}{\Gamma}& F' \ar@2[d]^{\sigma'}  \\
 H \ar@2[r]|-{\object@{|}}_{\gamma} & H'
} \quad \text{and} \quad \cd{
 G \ar@2[d]_{\tau} \ar@2[r]|-{\object@{|}}^{\beta} \dthreecell{dr}{\Delta}& G' \ar@2[d]^{\tau'}  \\
 K \ar@2[r]|-{\object@{|}}_{\delta} & K'
}
\]
in the hom-double categories $\mathfrak{Tricat}_3(\S, \T)$ and
$\mathfrak{Tricat}_3(\T, \U)$ respectively, we define their composite
\mbox{$\Delta \otimes \Gamma \colon \beta \otimes \alpha \Rrightarrow \delta
\otimes \gamma$} to be given by horizontally composing their underlying
families of cubical modifications in the locally cubical bicategory
$\mathfrak{Bicat}$:
\[(\Delta \otimes \Gamma)_{A, B} = \Gamma_{FA, FB} \otimes
\Delta_{A, B} \colon \beta_{FA, FB} \otimes \alpha_{A, B}
\Rrightarrow \delta_{FA, FB} \otimes \gamma_{A, B}\text.\] So in
particular, for any 1-cell $f \colon A \to B$ of $\S$, we have $(\Delta \otimes
\Gamma)_f$ given by the following pasting
\[\cd{
    GFf \ar@2[r]^{\beta_{Ff}} \ar@{=}[d] \dthreecell{dr}{\Delta_{Ff}} &
    G'Ff \ar@2[r]^{G'\alpha_{f}} \ar@{=}[d] \dthreecell{dr}{\tau'_{\alpha_f}} &
    G'F'f \ar@{=}[d] \\
    KFf \ar@2[r]^{\delta_{Ff}} \ar@{=}[d] \twoeq{dr} &
    K'Ff \ar@2[r]^{K'\alpha_{f}} \ar@{=}[d] \dthreecell{dr}{K'\Gamma_f} &
    K'F'f \ar@{=}[d] \\
    KHf \ar@2[r]_{\delta_{Hf}} &
    K'Hf \ar@2[r]_{K'\gamma_{f}} &
    K'H'f\text.}\]
Proving axioms (MA1) and (MA2) for this data amounts to constructing a further
succession of pasting equalities which traverse the interior of a $2 \times 2
\times 2$ cube, using:
\begin{itemize}
\item the corresponding axioms (MA1) or (MA2) for $\Delta$ and $\Gamma$,
\item the cubical modification axioms for the components of $\Delta$,
\item the icon axioms for the components of $\tau'$,
\item the pseudo-natural transformation axioms for the components of
    $\delta$,
\item and some further calculus with unnamed coherence isomorphisms.
\end{itemize}
Functoriality of this composition with respect to vertical composition is
inherited from that of horizontal composition of cubical modifications in
$\mathfrak{Bicat}$.

It remains to exhibit the pseudo-functoriality constraints for $\otimes$; so
let there be given lax homomorphisms and pseudo-icons
\[\cd{
  \S \ar@/^2em/[rr]^{F} \ar[rr]|{F'} \ar@/_2em/[rr]_{F''} & {}\dtwocell[0.36]{d}{\gamma} \dtwocell[-0.36]{d}{\alpha} &
  \T \ar@/^2em/[rr]^{G} \ar[rr]|{G'} \ar@/_2em/[rr]_{G''} & {}\dtwocell[0.36]{d}{\delta} \dtwocell[-0.36]{d}{\beta} &
  \U\text. \\
  & & & &
}\]\vskip-0.5\baselineskip\noindent We must exhibit invertible
globular icon modifications
\[i_{(G, F)} \colon \id_{GF} \Rrightarrow \id_G \otimes \id_F \colon GF \Rightarrow GF \quad \text{and} \quad
m_{(\beta, \alpha), (\delta, \gamma)} \colon (\delta \otimes \gamma)(\beta
\otimes \alpha) \Rrightarrow (\delta \beta) \otimes (\gamma \alpha)\text;\] and
to do this, we take their respective $(A, B)$-components to be the invertible
modifications witnessing pseudo-functoriality of horizontal composition in the
following diagram of homomorphisms and pseudo-natural transformations:
\[\cd{
  \S(A, B) \ar@/^2em/[rr]^{F_{A, B}} \ar[rr]|{F'_{A, B}} \ar@/_2em/[rr]_{F''_{A, B}} & {}\dtwocell[0.36]{d}{\gamma_{A, B}} \dtwocell[-0.36]{d}{\alpha_{A, B}} &
  \T(FA, FB) \ar@/^2em/[rr]^{G_{FA, FB}} \ar[rr]|{G'_{FA, FB}} \ar@/_2em/[rr]_{G''_{FA, FB}} & {}\dtwocell[0.36]{d}{\delta_{FA, FB}} \dtwocell[-0.36]{d}{\beta_{FA, FB}} &
  \sh{r}{0.75em} {\U(GFA, GFB)\text.} \\
  & & & &
}\] We must check that these data satisfy axioms (MA1) and (MA2). The proof is straightforward
manipulation using the pseudo-naturality axioms for $\delta$ and the
modification axioms for $\Pi^\delta$. Finally, the naturality of the maps
$m_{(\beta, \alpha), (\delta, \gamma)}$ in all variables follows componentwise;
as do the coherence axioms which $m$ and $i$ must satisfy.

In order to complete the definition of $\mathfrak{Tricat}_3$, all that remains  is to give the
associativity and unitality constraints for top-level composition, and to check the triangle
and pentagon axioms. At the level of 1-cells and vertical 2-cells, these are the corresponding
constraints from $\cat{Tricat}_2$; whilst at the level of horizontal 2-cells and 3-cells, we
suppose given trihomomorphisms and pseudo-icons
\[\cd{
  {\R} \ar@/^1em/[r]^{F} \ar@/_1em/[r]_{F'} \dtwocell{r}{\alpha} &
  {\S} \ar@/^1em/[r]^{G} \ar@/_1em/[r]_{G'} \dtwocell{r}{\beta} &
  {\T} \ar@/^1em/[r]^{H} \ar@/_1em/[r]_{H'} \dtwocell{r}{\gamma} &
  {\U\text,}
}\] and must exhibit an invertible pseudo-icon modification
\[
\cd[@+0.5em]{
 (HG)F \ar@2[d]_{a_{F, G, H}} \ar@2[r]|-{\object@{|}}^{(\gamma \otimes \beta) \otimes \alpha} \dthreecell{dr}{a_{\alpha, \beta, \gamma}}& (H'G')F' \ar@2[d]^{a_{F', G', H'}}  \\
 H(GF) \ar@2[r]|-{\object@{|}}_{\gamma \otimes (\beta \otimes \alpha)} &
 H'(G'F')
}\] where $a_{F, G, H}$ and $a_{F', G', H'}$ are the corresponding constraints
from $\cat{Tricat}_2$. So we take the $(A, B)$th component of this pseudo-icon
modification to be the cubical modification providing the associativity
constraint for the composition
\[\cd[@C+2.3em]{
  {\R(A, B)} \ar@/^1em/[r]^{F_{A, B}} \ar@/_1em/[r]_{F'_{A, B}} \dtwocell[0.45]{r}{\alpha_{A, B}} &
  {\S(FA, FB)} \ar@/^1em/[r]^{G_{FA, FB}} \ar@/_1em/[r]_{G'_{FA, FB}} \dtwocell[0.37]{r}{\beta_{FA, FB}} &
  {\T(GFA, GFB)} \ar@/^1em/[r]^{H_{GFA, GFB}} \ar@/_1em/[r]_{H'_{GFA, GFB}} \dtwocell[0.32]{r}{\gamma_{GFA, GFB}} &
  {\U(HGFA, HGFB)}
}\] in the locally cubical bicategory $\mathfrak{Bicat}$. We must check that
these data satisfy the axioms for an icon modification; let us do only (MA2),
since (MA1) is similar. We first observe that the 3-cells $\Pi^{(\gamma \otimes
\beta) \otimes \alpha}_A$ and $\Pi^{\gamma \otimes (\beta \otimes \alpha)}_A$
are obtained by pasting together what is essentially the same $3 \times 3$
diagram of 3-cells, and some trivial calculus with unnamed coherence cells
shows that they are precisely the same diagram, modulo rewriting of the
boundary, so that the latter 3-cell may be obtained from the former by pasting
with unnamed coherence cells. But this is precisely the content of axiom (MA2).

Finally, we must check that these icon modifications $a_{\alpha, \beta, \gamma}$ are natural
in in $\alpha$, $\beta$ and $\gamma$, and satisfy the pentagon and triangle equalities.  Each
of these follows componentwise from the corresponding facts in $\mathfrak{Bicat}$. This
completes the definition of $\mathfrak{Tricat}_3$.

\section{From locally cubical bicategories to tricategories}\label{Sec:fromltt}
In the previous Section, we constructed a locally cubical bicategory of
tricategories which we called $\mathfrak{Tricat}_3$. Recall from
Section~\ref{Sec:towards} that one reason for doing this was so that we could
deduce the existence of the tricategory $\cat{Tricat}_3$. The purpose of this
section is to describe the general machinery which will allow us to do this.

The construction takes a well-behaved locally cubical bicategory $\mathfrak B$ and builds a
tricategory out of it. This tricategory will have the same 0-\ and 1-cells as $\mathfrak B$;
as 2-cells, the horizontal 2-cells of $\mathfrak B$; and as 3-cells, the globular 3-cells of
$\mathfrak B$. The main point of interest is the construction of the tricategorical
associativity constraints, which are to be given by \emph{horizontal} 2-cells of $\mathfrak
B$. Since the associativity constraints in $\mathfrak B$ are given by \emph{vertical} 2-cells,
we will need some kind of linkage between the two types of 2-cell in order to proceed.

\begin{Defn}
A pseudo double category $\mathfrak C$ is \defn{fibrant} if the functor $(s, t)
\colon \C_1 \to \C_0 \times \C_0$ is an isofibration.
\end{Defn}
\noindent Recall here that a functor $F \colon \A \to \B$ between categories is
an \emph{isofibration} if whenever we have a object $a \in \A$ and isomorphism
$\phi \colon Fa \to b$ in $\B$, there exists an object $c \in \A$ and
isomorphism $\theta \colon a \to c$ such that $Fc = b$ and $F\theta = \phi$.
Thus a pseudo double category $\mathfrak C$ is fibrant just when every diagram
like (a) below with $f$ and $g$ isomorphisms has a filler like (b) for which the 2-cell $\theta$ is
invertible as an arrow of $\C_1$:
\[\text{(a)} \quad \cd{
x \ar[d]_{f} & x' \ar[d]^{g}  \\
y \ar[r]|-{\object@{|}}_{k} & y'} \qquad \rightsquigarrow \qquad \text{(b)} \quad \cd{
x \ar[d]_{f} \ar[r]|-{\object@{|}}^{h} \dtwocell{dr}{\theta}& x' \ar[d]^{g}  \\
y \ar[r]|-{\object@{|}}_{k} & y'\text.}\] Thus fibrancy is precisely the
property which \cite{GP1} refers to as \emph{horizontal invariance}. We may
reformulate this property in various useful ways, and since detailed accounts
of this process may be found in \cite{tomf} or \cite{GP2}, we restrict
ourselves here to recording those equivalent formulations which will be useful
to us.

For the first, we consider the case of the above filling condition where $g$
and $k$ are both identities: given a vertical map $f \colon x \to y$ of
$\mathfrak C$, it asserts the existence of a horizontal 1-cell $\overline f$
and a 2-cell $\epsilon_f$ fitting into the diagram:
\[\cd{
x \ar[d]_{f} \ar[r]|-{\object@{|}}^{\overline f} \dtwocell{dr}{\epsilon_f}& y \ar[d]^{\id_y}  \\
y \ar[r]|-{\object@{|}}_{\DI y} & y\text.}\] From this, we may define a further 2-cell $\eta_f$ as
the composite
\[\cd{
  x \ar[d]_{\id_x} \ar[r]|-{\object@{|}}^{\DI x} \dtwocell{dr}{\eta_f} &
  x \ar[d]^{f} \\
  x \ar[r]|-{\object@{|}}_{\overline f} & y} \quad := \quad \cd{
  x \ar[d]_{f} \ar[r]|-{\object@{|}}^{\DI x} \dtwocell{dr}{\DI f} &
  x \ar[d]^{f} \\
  y \ar[d]_{f^{-1}} \ar[r]|-{\object@{|}}^{\overline f} \dtwocell{dr}{\epsilon_f^{-1}} &
  y \ar[d]^{\id_y} \\
  x \ar[r]|-{\object@{|}}_{\overline f} & y\text.}\]
Now the pair $(\eta_f, \epsilon_f)$ satisfy the triangle identities:
\[\epsilon_f . \eta_f = \DI f \colon \DI x \Rightarrow \DI y \quad \text{and} \quad \epsilon_f \ast \eta_f = (l^{-1}r)_{\overline f} \colon \overline f . \DI x \Rightarrow \DI y . \overline f\text,\]
and so, in the terminology of \cite{GP2}, $f$ and $\overline f$ are
\emph{orthogonal companions}; which gives us the ``only if'' direction of:
\begin{Prop}
A pseudo double category $\mathfrak C$ is fibrant iff every vertical
isomorphism has an orthogonal companion.
\end{Prop}
\noindent For the ``if'' direction, suppose that we are given a diagram like
(a); then we can complete it to a diagram like (b) by taking $h$ to be
$\overline{g^{-1}} . (k . \overline f)$ and $\theta$ to be the 2-cell:
\[\overline{g^{-1}} . (k . \overline f) \xRightarrow{(\DI g .\epsilon_{g^{-1}}) \ast (\id_k \ast
\epsilon_f)} \DI{y'} . (k . \DI y) \xRightarrow{\DI{y'} \ast r_k} \DI{y'} . k \xRightarrow{l_k}
k\text.\]

Thus each of $\mathfrak{Cat}$, $\mathfrak{Rng}$ and $\mathfrak{Span}(\C)$ is a
fibrant double category: for $\mathfrak{Cat}$, the horizontal companion of a
functor $F \colon \C \to D$ is the profunctor $\overline F(\mathord \thg,
\mathord ?) = \D(\mathord \thg, F\mathord ?)$; for $\mathfrak{Rng}$, the
companion of a homomorphism $f \colon R \to S$ is $S$ itself, viewed as a left
$S$-, right $R$-module; and in $\mathfrak{Span}(\C)$, the companion of a
morphism $f \colon C \to D$ is the span $C \stackrel{id}{\leftarrow} C
\stackrel{f}{\rightarrow} D$. Observe that in all of these examples, it is
arbitrary vertical morphisms, and not just the isomorphisms, which have
companions: such pseudo double categories are essentially the \emph{pro-arrow
equipments} of \cite{wood1, wood2}. A more detailed analysis of this
correspondence may be found in Appendix C of \cite{mike}.

\begin{Prop}
Let $\mathfrak C$ be a fibrant double category equipped with a choice of
orthogonal companion for every vertical isomorphism. Then the assignation $f
\mapsto \overline f$ underlies an identity-on-objects homomorphism of
bicategories
\[(\overline{\ \ \mathstrut})\colon V^{\text{iso}}(\mathfrak C) \to H(\mathfrak C)\text,\]
where $V^{\text{iso}}(\mathfrak C)$ is the category of objects and vertical
isomorphisms in $\mathfrak C$. Moreover, if we are given vertical isomorphisms
$f \colon w \to y$ and $g \colon x \to z$ in $\mathfrak C$, then pasting with
$\eta_f$ and $\epsilon_g$ induces a bijection between the set of 2-cells of the
form (c) and the set of 2-cells of the form (d):
\[(c) \quad \cd{
x \ar[d]_{f} \ar[r]|-{\object@{|}}^{h} \dtwocell{dr}{\alpha} & x' \ar[d]^{g}  \\
y \ar[r]|-{\object@{|}}_{k} & y'} \qquad \text{and} \qquad \text{(d)} \quad \cd{
x \ar[d]_{\id_x} \ar[r]|-{\object@{|}}^{\overline g . h} \dtwocell{dr}{\overline \alpha} & y' \ar[d]^{\id_{y'}}  \\
x \ar[r]|-{\object@{|}}_{k . \overline f} & y'\text;}\] and $\overline \alpha$ is invertible
as an arrow of $\C_1$ just when $\alpha$ is. Furthermore, these bijections satisfy four
evident axioms expressing their functoriality with respect to vertical and horizontal
composition of 2-cells.\end{Prop} If we remove the restriction to vertical isomorphisms, then
the structure described in this Proposition is that of a \emph{pseudo folding structure} in
the sense of \cite{tomf}. The proof of the Proposition is straightforward manipulation, and it
is not hard to prove a converse~--~namely, that from a homomorphism of bicategories
$(\overline{\ \ \mathstrut})\colon V^{\text{iso}}(\mathfrak C) \to H(\mathfrak C)$ and a
bijective assignation $\alpha \mapsto \overline \alpha$ on 2-cells satisfying the four
functoriality axioms, one may define a choice of orthogonal companion for every vertical
isomorphism. A proof of this correspondence may be extracted from the pages leading up to
Theorem 3.28 of \cite{tomf}.

\begin{Defn}
The 2-category $\cat{DblCat}_f$ has objects being fibrant pseudo double
categories equipped with a choice of orthogonal companions; as 1-cells, the
homomorphisms between the underlying double categories; and as 2-cells, the
vertical transformations between them.
\end{Defn}
\noindent One may reasonably ask why we do not require the 1-cells $F \colon
\mathfrak C \to \mathfrak D$ of $\cat{DblCat}_f$ to respect the choices of
orthogonal companions in $\mathfrak C$ and $\mathfrak D$. The reason is that,
in fact, \emph{any} homomorphism between objects of $\cat{DblCat}_f$ will
automatically respect these choices in a unique way. To make this explicit, let
us say that a homomorphism $F \colon \mathfrak C \to \mathfrak D$ between
objects of $\cat{DblCat}_f$ is a
\defn{fibrant homomorphism} if, for every invertible vertical 1-cell $f \colon x \to y$ of $\mathfrak C$, there is given
an invertible globular 2-cell \[\mu_f \colon F(\overline f) \Rightarrow
\overline{Ff} \colon Fx \tor Fy\] of $\mathfrak D$, subject to three axioms.
The first two equate, respectively, the two possible 2-cells in $\mathfrak D$
from $\DI {Fx}$ to $\overline{F(\id_x)}$; and from $F(\overline g) \, . \, F
(\overline f)$ to $\overline{F(gf)}$. The third axiom concerns a 2-cell
$\alpha$ of the type (c) above, and equates the two globular 2-cells
\[
\cd[@-1em]{
Fx \ar@/_4em/[dd]|-{\object@{|}}_{\overline {Ff}} \ar[dd]|-{\object@{|}}_{F(\overline f)} \ar[rr]|-{\object@{|}}^{Fh} \dtwocell{ddrr}{F(\overline \alpha)} && x' \ar[dd]|-{\object@{|}}^{F(\overline g)}  \\
{}\ltwocell[-0.9]{r}{\mu_f} &
\\
Fy \ar[rr]|-{\object@{|}}_{Fk} && Fy'} \quad \text{and} \quad \cd[@-1em]{
Fx \ar[dd]|-{\object@{|}}_{\overline {Ff}} \ar[rr]|-{\object@{|}}^{Fh} \dtwocell{ddrr}{\overline {F\alpha}} && x' \ar[dd]|-{\object@{|}}^{\overline {Fg}} \ar@/^4em/[dd]|-{\object@{|}}^{F(\overline g)}  \\
& & {}\ltwocell[-0.9]{l}{\mu_f}
\\
Fy \ar[rr]|-{\object@{|}}_{Fk} && Fy'\text.}
\]
Now, given any homomorphism $F \colon \mathfrak C \to \mathfrak D$ between
objects of $\cat{DblCat}_f$, we may make it into a fibrant homomorphism as
follows. Given an invertible vertical arrow $f \colon x \to y$ of $\mathfrak
C$, we can consider the globular 2-cell
\[\overline{F\epsilon_f} \colon \overline{\id_{Fy}} . F (\overline f)
\Rightarrow F\DI y . \overline{Ff}\] of $\mathfrak D$; and since both $\overline{\id_{Fy}}$
and $F\DI y$ are isomorphic to $\DI {Fy}$, we obtain from this a globular 2-cell $\mu_f \colon
F (\overline f) \Rightarrow \overline{Ff}$, which is easily checked to satisfy the three
axioms. And in fact, this is the only possible structure of fibrant homomorphism on $F$: for
given an arbitrary such structure, applying the third axiom to the 2-cells $\epsilon_f$ in
$\mathfrak C$ shows that the maps $\mu_f$ must coincide with those defined above. A
similar argument applies to the 2-cells of $\cat{DblCat}_f$.

A conceptual explanation of why this should be the case is that
$\cat{DblCat}_f$ is, in some sense, the 2-category of algebras for a
particularly simple kind of 2-dimensional monad on $\cat{DblCat}$, the kind
which \cite{propertylike} calls \emph{pseudo-idempotent}: and such monads have
the property that the forgetful functor from the 2-category of algebras and
algebra pseudomorphisms to the underlying base 2-category is 2-fully faithful.
The qualifier ``in some sense'' covers a slight wrinkle in this explanation:
namely, that the 2-monad which gives rise to $\cat{DblCat}_f$ lives not on
$\cat{DblCat}$ but on $\cat{DblCat}_\text{str}$, the 2-category of pseudo
double categories and \emph{strict} homomorphisms between them, so that making
this argument rigourous would require a little more work.

\begin{Defn}
We will say that a locally cubical bicategory is \defn{locally fibrant} just
when each of its hom-double categories is fibrant.
\end{Defn}
In particular, a monoidal double category is locally fibrant just when its
underlying pseudo double category is fibrant, so that all of our examples of
monoidal double categories are locally fibrant. The locally cubical bicategory
$\mathfrak{DblCat}$ is easily seen \emph{not} to be locally fibrant; on the
other hand, we may show that, for pseudo double categories $\mathfrak C$ and
$\mathfrak D$, if $\mathfrak D$ is fibrant then so is $\HOM{\mathfrak C,
\mathfrak D}$. It follows that the locally cubical bicategory
$\mathfrak{DblCat}_f$, of fibrant double categories and all cells between them,
is itself locally fibrant; and since any bicategory is trivially fibrant, that
the locally cubical bicategory $\mathfrak{Bicat}$ is too.

We will now show that every locally fibrant locally cubical bicategory gives
rise to a tricategory. We begin with a technical result:

\begin{Prop}\label{trihomom}
Let  $\cat{DblCat}_g$ be the maximal sub-2-category of $\cat{DblCat}_f$ with
only invertible 2-cells. Then the functor of mere categories $\cat{DblCat}_g
\hookrightarrow \cat{DblCat} \xrightarrow \HH \cat{Bicat}$ can be extended to a
trihomomorphism
\[\HH \colon \cat{DblCat}_g \to \cat{Bicat}\text.\]
\end{Prop}
\begin{proof}
First we define $\HH$ on cells. This is already done for 0-\ and 1-cells, and
since $\cat{DblCat}_g$ has no non-trivial 3-cells, it remains only to define it
on 2-cells. So let there be given an invertible vertical transformation $\alpha
\colon F \Rightarrow G \colon \mathfrak C \to \mathfrak D$. We define a
pseudo-natural transformation $\HH\alpha \colon \HH F \Rightarrow \HH G$ by
taking
\[(\HH \alpha)_x = \overline{\alpha_x} \colon Fx \to Gx\text, \qquad \text{and} \qquad (\HH \alpha)_f = \cd{
Fx \ar[d]_{\overline{\alpha_x}} \ar[r]^{Ff} \dtwocell{dr}{\overline{\alpha_f}}& Fy \ar[d]^{\overline{\alpha_y}}  \\
Gx \ar[r]_{Gf} & Gy\text.}\] The transformation axioms for $\HH\alpha$ follow straightforwardly
from the vertical transformation axioms for $\alpha$ and the functoriality of
$(\overline{\ \ \mathstrut})$ with respect to 2-cell composition. Next we
ensure that $\HH$ is locally a homomorphism of bicategories, which entails
giving modifications $i_F \colon \id_{\HH F} \Rrightarrow \HH(\id_F)$ and
$m_{\alpha, \beta} \colon
 \HH \beta . \HH \alpha\Rrightarrow \HH(\beta . \alpha)$. These will have 2-cell components
\[(i_F)_x \colon \id_{Fx} \Rightarrow \overline{\id_{Fx}} \qquad \text{and} \qquad (m_{\alpha, \beta})_x \colon \overline{\beta_x} \, . \, \overline{\alpha_x} \Rightarrow \overline{\beta_x \alpha_x}\]
in $\HH \mathfrak D$ given by the pseudo-functoriality constraints for
$(\overline{\ \ \mathstrut})$. The coherence axioms for these data therefore
follow pointwise. Next, we must give adjoint pseudo-natural equivalences
\[{\b \chi}_{\mathfrak C, \mathfrak D, \mathfrak E} \colon \HH(\thg) \otimes \HH(?) \Rightarrow \HH(\thg \otimes ?) \colon \cat{DblCat}_g(\mathfrak C, \mathfrak
D) \times \cat{DblCat}_g(\mathfrak B, \mathfrak C) \to \cat{Bicat}(\HH \mathfrak B, \HH \mathfrak
D)\text.\] Observe that the homomorphisms $\HH(\thg) \otimes \HH(?)$ and $\HH(\mathord \thg \otimes
\mathord ?)$ agree on objects, and thus we may consider icons between them: in
particular, any \emph{invertible} icon between them will give rise to an
adjoint pseudo-natural equivalence, and so to give $\chi$ it suffices to give
invertible icons $\chi \colon \HH(\thg) \otimes \HH(?) \Rightarrow \HH(\mathord
\thg \otimes \mathord ?)$. So consider a pair of horizontally composable
2-cells
\[\cd{
  {\mathfrak B} \ar@/^1em/[r]^{F} \ar@/_1em/[r]_{F'} \dtwocell{r}{\alpha} &
  {\mathfrak C} \ar@/^1em/[r]^{G} \ar@/_1em/[r]_{G'} \dtwocell{r}{\beta} &
  {\mathfrak D}
}\]in $\cat{DblCat}_g$: we must give a modification $\chi_{\alpha, \beta} \colon \HH \beta \ast \HH
\alpha \Rrightarrow \HH(\beta \ast \alpha)$. Now, these two pseudo-natural
transformations have respective $x$-components given by \begin{align*} (\HH
\beta \ast \HH \alpha)_x & = GFx \xrightarrow{\overline{\beta_{Fx}}} G'Fx
\xrightarrow{G'\overline{\alpha_x}} G'F'x\\ \text{and\ \ } \HH(\beta \ast
\alpha)_x & = GFx \xrightarrow{\overline{G'{\alpha_x} . \beta_{Fx}}}
G'F'x\text,\end{align*} and so we take $(\chi_{\alpha, \beta})_x$ to be the
2-cell
\[G'\overline{\alpha_x} \, . \, \overline{\beta_{Fx}} \xRightarrow{\mu_{\alpha_x} . 1} \overline{G'\alpha_x} \, . \, \overline{\beta_{Fx}} \xRightarrow{\ \cong\ }  \overline{G'\alpha_x . \beta_{Fx}}\]
of $\HH \mathfrak D$. The modification axioms for $\chi_{\alpha, \beta}$ follow
from the third fibrant homomorphism axiom and the functoriality axioms for
$(\overline{\ \ \mathstrut})$ with respect to 2-cell composition. We must
verify that these components $\chi_{\alpha, \beta}$ satisfy the three axioms
making $\chi$ into an icon. The first is vacuous, whilst the second and third
follow by a diagram chase using the axioms for a fibrant homomorphism. We argue
entirely analogously in order to give the adjoint equivalences $\iota \colon
I_{\HH x} \Rightarrow \HH(I_x)$.

Next we must give invertible modifications $\omega$, $\delta$ and $\gamma$. In
the case of $\omega$, for instance, this involves giving invertible
modifications
\[\cd[@!C@!R@R+1em]{
    & (\HH(\thg)\otimes\HH(?))\otimes\HH(\mathord \ast)
      \ar@2[dl]_{\chi \otimes 1} \ar@2[dr]^{\a} \\
    \HH(\mathord \thg \otimes \mathord ?)\otimes\HH(\mathord\ast)
      \ar@2[d]_{\chi} \rthreecell{drr}{\omega} & &
    \HH(\thg) \otimes (\HH(?) \otimes \HH(\mathord \ast))
      \ar@2[d]^{1 \otimes \chi} \\
    \HH((\mathord \thg \otimes \mathord ?) \otimes \mathord \ast)
      \ar@2[dr]_{\HH\a} & &
    \HH(\thg) \otimes \HH(\mathord ? \otimes \mathord \ast)
      \ar@2[dl]^{\chi} \\ &
    \HH(\thg \otimes (\mathord ? \otimes \mathord \ast))\text.
}\] To do this, observe first that every pseudo-natural transformation bounding this diagram may
also be viewed as an icon. We already know this for $\chi$ and hence also for
$1 \otimes \chi$ and $\chi \otimes 1$; and it is so for $\a$ and $\HH \a$ since
composition of 1-cells in both $\cat{DblCat}_g$ and $\cat{Bicat}$ is
\emph{strictly} associative. If we now compose all the 2-arrows in this diagram
\emph{qua} icons, we obtain two further icons $\sigma, \tau \colon
(\HH(\thg)\otimes\HH(?))\otimes\HH(\mathord \ast) \Rightarrow \HH(\thg \otimes
(\mathord ? \otimes \mathord \ast))$: and a long but straightforward diagram
chase with the fibrant homomorphism axioms shows that these two icons are, in
fact, equal.

On the other hand, if we compose the two sides \emph{qua} pseudo-natural
transformations, then the pseudo-naturals that we get will not necessarily be
icons, but they will, at least, be \emph{isomorphic} to icons, namely the icons
$\sigma$ and $\tau$ respectively. Thus we take $\omega$ to be the composite of
the invertible modification from the left-hand side of this diagram to $\sigma
= \tau$ and the invertible modification from $\tau$ to the right-hand side. We
proceed similarly for the invertible modifications $\delta$ and $\gamma$.

The final thing to check are the two trihomomorphism axioms, equating certain
pastings of 3-cells in $\cat{Bicat}$. But all the 3-cells in question are
either coherence 3-cells of $\cat{Bicat}$; or component 3-cells of $\omega$,
$\delta$ and $\gamma$. But these latter 3-cells are in turn built from
coherence 3-cells of $\cat{Bicat}$ and coherence 3-cells for the local
homomorphisms $(\overline{\ \ \mathstrut})$. The result thus follows by
coherence for tricategories and bicategorical coherence for functors.
\end{proof}
\begin{Thm}\label{locallydoubletotri}
Let $\mathfrak C$ be a locally fibrant locally cubical bicategory with chosen
companions in each hom. Then there is a tricategory $\T$ with the same objects
as $\mathfrak C$, and
\[\T(A, B) = \HH\big(\mathfrak C(A, B)\big)\text.\]
\end{Thm}
\begin{proof}
We begin by observing that both $\cat{DblCat}_g$ and $\cat{Bicat}$ come
equipped with finite product structure; and that the trihomomorphism $\HH$
preserves the cartesian product of $j$-cells for $j = 0, 1, 2, 3$. Now, the
top-level composition and identity functors for $\T$ are given by applying
$\HH$ to the corresponding data (LDD3) and (LDD4) for $\mathfrak C$:
\[1 = \HH 1 \xrightarrow{\HH \elt{I_A}} \HH\big(\mathfrak C(A, A)) = \T(A, A)\]
and
\[\T(B, C) \times \T(A, B) = \HH\big(\mathfrak C(B, C) \times \mathfrak C(A, B)\big) \xrightarrow{\HH \mathord \otimes} \HH\big(\mathfrak C(A, C)\big) = \T(A, C)\text.\]
To obtain the pseudo-natural adjoint equivalences $\a$, $\l$ and $\r$
witnessing the associativity and unitality of this composition, we apply $\HH$
to the corresponding data (LDD5) and (LDD6) for $\mathfrak C$. Since each of
$a$, $l$ and $r$ is an adjoint equivalence (in fact, an isomorphism) in the
relevant hom of $\cat{DblCat}_g$, the same will obtain for their images in
$\cat{Bicat}$; and because $\HH$ strictly preserves both cartesian products and
composition of 1-cells, these adjoint equivalences will have the correct
sources and targets.

Next we must give the invertible modifications $\pi$, $\mu$, $\lambda$ and
$\rho$. To obtain $\pi$, for example, we begin by applying $\HH$ to the axiom
(LDA2) for $\mathfrak C$. This yields an equality of 2-cells in $\cat{Bicat}$;
however, these 2-cells are not of the right form to be the source and target of
$\pi$. In order to make them so, we may adjust by coherence 3-cells in
$\cat{Bicat}$ whose existence is guaranteed by the coherence theorem for
trihomomorphisms. Consequently, we may take $\pi$ to be given by the composite
of these coherence 3-cells; and similarly for $\mu$, $\lambda$ and $\rho$.

Finally, we must check the three tricategory axioms. These are normally stated
in a ``local'' form, asserting the equality of certain pastings of 3-cells in
the relevant hom-bicategories; but in this situation, it will be more
appropriate to consider them in their ``global'' form. Each such axiom amounts
to giving a diagram of 2-\ and 3-cells in $\cat{Bicat}$, whose vertices are
pasting diagrams built from copies of the 2-cells $a$, $l$ and $r$, and whose
arrows are 3-cells between those 2-cells, built from copies of $\pi$, $\mu$,
$\lambda$ and $\rho$; and asserting that the two ways around this diagram
coincide.

To show this, we consider the corresponding diagram for $\mathfrak C$. This is
a diagram of 2-\ and 3-cells in $\cat{DblCat_g}$, which since $\cat{DblCat}_g$
has only identity 3-cells, must commute. Hence by applying $\HH$ we obtain a
commutative diagram in $\cat{Bicat}$, which, unfortunately, has both the wrong
vertices and the wrong arrows. Nonetheless, by the coherence theorem for
functors, each ``wrong'' vertex admits an isomorphism 3-cell to the ``right''
vertex; and in such a way that composing these isomorphism 3-cells with the
``wrong'' arrows yields the ``right'' arrows.
\end{proof}

Special cases of this theorem give us new proofs of some existing results.
Restricting to the one-object case, we have the result that \emph{any fibrant
monoidal double category gives rise to a monoidal bicategory}; this statement
and a sketch proof appear as Theorem B.4 of \cite{mike}. In particular, we
obtain elegant proofs that the bicategories of rings and bimodules, of
categories and profunctors, and of spans internal to a cartesian category $\C$
are all monoidal bicategories.\footnote{For the last two of these examples, the
machinery of \cite{cb2} provides another elegant proof of this fact, and in
fact goes further, showing that the monoidal bicategories in question are
\emph{symmetric} monoidal bicategories.}\ \  Finally, applying this theorem to
the fibrant locally cubical bicategory $\mathfrak{Bicat}$, we deduce the
existence of a tricategory of bicategories $\cat{Bicat}$. Again, the result is
not new, but the proof is, showing how the tricategory structure on
$\cat{Bicat}$ may be induced from a piece of canonical, universally determined
structure, namely the biclosed structure on $\cat{DblCat}$.

\section{A tricategory of tricategories}\label{tritri}
We are now finally in a position to prove Theorem \ref{mainresultprop}, which
asserts the existence of the tricategory of tricategories $\cat{Tricat}_3$. We
will do this by applying the machinery of the previous section to the locally
cubical bicategory $\mathfrak{Tricat}_3$. In order to do this, we must first
prove that $\mathfrak{Tricat}_3$ is locally fibrant.
\begin{Prop}
Each pseudo double category $\mathfrak{Tricat}_3(\S, \T)$ is fibrant.
\end{Prop}
\begin{proof}
Suppose we are given an invertible ico-icon $\alpha \colon F \Rightarrow G
\colon \S \to \T$. We must provide a pseudo-icon $\overline \alpha \colon F
\Tor G$ and an invertible icon modification
\[\cd{
 F \ar@2[d]_{\alpha} \ar@2[r]|-{\object@{|}}^{\overline \alpha} \dthreecell{dr}{\epsilon_\alpha}& G \ar@2[d]^{\id_G}  \\
 G \ar@2[r]|-{\object@{|}}_{\DI G} & G.
}\] Now by Proposition \ref{embed}, there is a bijection between the ico-icons $F \Rightarrow G$
and the oplax icons $F \Tor G$ with identity 2-cell components: for which the
invertible ico-icons on the one side correspond to the pseudo-icons on the
other. Thus we take $\overline \alpha$ to be the pseudo-icon corresponding to
$\alpha$ under this bijection.

To give the icon modification $\epsilon_\alpha$, we must give 3-cells
$(\epsilon_\alpha)_f \colon {\overline \alpha}_f \Rrightarrow (\DI G)_f$ of
$\T$, forming the components of an $\ob \S \times \ob \S$-indexed family of
cubical modifications, and satisfying axioms (MA1) and (MA2). Since we have
${\overline \alpha}_f = (\DI G)_f = \id_{Ff} \colon Ff \Rightarrow Ff$, we take
$(\epsilon_\alpha)_f = \id_{\id_f}$. The cubical modification axioms and axioms
(MA1) and (MA2) now follow by coherence for bicategories.
\end{proof}
And so finally we obtain:
\begin{Cor}\label{finalresult}
There is a tricategory $\cat{Tricat}_3$ with objects being tricategories;
1-cells, lax homomorphisms; 2-cells, pseudo-icons; and 3-cells, pseudo-icon
modifications.
\end{Cor}
\begin{Cor}
The tricategory $\cat{MonBicat}$ of monoidal bicategories, weak monoidal
functors, weak monoidal transformations, and monoidal modifications is
triequivalent to the full sub-tricategory of $\cat{Tricat}_{3}$ consisting of
those tricategories with a single object.
\end{Cor}


\begin{thebibliography}{10}

\bibitem{bicats}
{\sc B{\'e}nabou, J.}
\newblock Introduction to bicategories.
\newblock In {\em Reports of the Midwest Category Seminar}. Springer, Berlin,
  1967, pp.~1--77.

\bibitem{cb2}
{\sc Carboni, A., Kelly, G.~M., Walters, R., and Wood, R.}
\newblock Cartesian bicategories {II}.
\newblock {\em Theory and Applications of Categories 19\/} (2008), 
  93--124.

\bibitem{Carmody}
{\sc Carmody, S.}
\newblock {\em Cobordism Categories}.
\newblock PhD thesis, Cambridge, 1995.

\bibitem{cheng-gurski} 
{\sc Cheng, E., and Gurski, N.}
\newblock The periodic table of n-categories II: degenerate tricategories.
\newblock {\em Cahiers de Topologie et Geom\'etrie Diff\'erentielle Cat\'egoriques}, in press.

\bibitem{Biclosed}
{\sc Day, B.}
\newblock Biclosed bicategories: localisation of convolution.
\newblock Unpublished, available as arXiv preprint 0705.3485.

\bibitem{monoidalbicats}
{\sc Day, B., and Street, R.}
\newblock Monoidal bicategories and {H}opf algebroids.
\newblock {\em Advances in Mathematics 129}, 1 (July 1997), 99--157.

\bibitem{tomf}
{\sc Fiore, T.}
\newblock Pseudo algebras and pseudo double categories.
\newblock {\em Journal of Homotopy and Related Structures 2\/} (2007), 119--170.

\bibitem{doubleclubs}
{\sc Garner, R.}
\newblock Double clubs.
\newblock {\em Cahiers de Topologie et Geom{\'e}trie Diff{\'e}rentielle
  Cat{\'e}goriques 47}, 4 (2006), 261--317.

\bibitem{GPS}
{\sc Gordon, R., Power, A.~J., and Street, R.}
\newblock Coherence for tricategories.
\newblock {\em Memoirs of the American Mathematical Society 117}, 558 (1995).

\bibitem{GP1}
{\sc Grandis, M., and Par{\'e}, R.}
\newblock Limits in double categories.
\newblock {\em Cahiers de Topologie et Geom{\'e}trie Diff{\'e}rentielle
  Cat{\'e}goriques 40}, 3 (1999), 162--220.

\bibitem{GP2}
{\sc Grandis, M., and Par{\'e}, R.}
\newblock Adjoints for double categories.
\newblock {\em Cahiers de Topologie et Geom{\'e}trie Diff{\'e}rentielle
  Cat{\'e}goriques 45}, 3 (2004), 193--240.

\bibitem{nicktricats}
{\sc Gurski, N.}
\newblock {\em An algebraic theory of tricategories}.
\newblock PhD thesis, University of Chicago, 2006.

\bibitem{propertylike}
{\sc Kelly, G.~M., and Lack, S.}
\newblock On property-like structures.
\newblock {\em Theory and Applications of Categories 3\/} (1997), No. 9,
  213--250.

\bibitem{review}
{\sc Kelly, G.~M., and Street, R.}
\newblock Review of the elements of {$2$}-categories.
\newblock In {\em Category Seminar (Proc. Sem., Sydney, 1972/1973)}. Springer,
  Berlin, 1974, pp.~75--103. Lecture Notes in Math., Vol. 420.

\bibitem{Lackicons}
{\sc Lack, S.}
\newblock Icons.
\newblock {\em Applied Categorical Structures 18\/} (2010), No. 3,
  289--307.
  
  \bibitem{Lack}
{\sc Lack, S.}
\newblock {\em The algebra of distributive and extensive categories}.
\newblock PhD thesis, Cambridge, 1995.



\bibitem{LP}
{\sc Lack, S., and Paoli, S.}
\newblock 2-nerves for bicategories.
\newblock {\em K Theory 38}, 2 (2008), 153--175.

\bibitem{MLP}
{\sc Mac~Lane, S., and Par{\'e}, R.}
\newblock Coherence for bicategories and indexed categories.
\newblock {\em Journal of Pure and Applied Algebra 37}, 1 (1985), 59--80.

\bibitem{mike}
{\sc Shulman, M.}
\newblock Framed bicategories and monoidal fibrations.
\newblock {\em Theory and Application of Categories 20\/} (2008), 650--738.
  
\bibitem{Enrichedbicats}
{\sc Street, R.}
\newblock Enriched categories and cohomology.
\newblock {\em Reprints in Theory and Applications of Categories}, 14 (2005),
  1--18.
\newblock Reprinted from Quaestiones Math. {\bf 6} (1983), no. 1-3, 265--283,
  with new commentary by the author.

\bibitem{streetomega} 
{\sc Street, R.}
\newblock Weak omega-categories.
\newblock ``Diagrammatic Morphisms and Applications,'' {\em Contemporary Mathematics},
 318 (2003), 207--213.


\bibitem{dv} 
{\sc Verity, D.}
\newblock {\em Enriched categories, internal categories, and change of base}.
\newblock PhD thesis, Cambridge, 1992.

\bibitem{wood1}
{\sc Wood, R.~J.}
\newblock Abstract proarrows {I}.
\newblock {\em Cahiers de Topologie et Geom{\'e}trie Diff{\'e}rentielle
  Cat{\'e}goriques 23}, 3 (1982), 279--290.

\bibitem{wood2}
{\sc Wood, R.~J.}
\newblock Proarrows {II}.
\newblock {\em Cahiers de Topologie et Geom{\'e}trie Diff{\'e}rentielle
  Cat{\'e}goriques 26}, 2 (1985), 135--168.
  
\end{thebibliography}
\end{document}